\documentclass[11pt, oneside]{amsart}
\usepackage{geometry} % see geometry.pdf on how to lay out the page. There's lots.
\usepackage{graphicx}
\usepackage{amsmath}
\usepackage{amsfonts}
\usepackage{mathrsfs}
\usepackage{amssymb}
\usepackage{wrapfig}
\usepackage{bbold}
\usepackage{xcolor, soul}
\usepackage{mdframed}
\usepackage{enumitem}
\usepackage[colorlinks=true, pdfstartview=FitV, linkcolor=blue, 
        citecolor=blue, urlcolor=blue]{hyperref}

\newtheorem{thm}{Theorem}
\newtheorem{lem}{Lemma}[section]
\newtheorem{cor}[lem]{Corollary}
\newtheorem{prop}[lem]{Proposition}

\theoremstyle{remark}
\newtheorem{remark}[lem]{Remark}
\AtEndEnvironment{remark}{\null\hfill\rule{.5ex}{1.4ex}\,}%\square

\newenvironment{pf}{\noindent {\em Proof}.\ \ }{\hspace*{\fill}\rule{.5ex}{1.4ex}\,}

\numberwithin{equation}{section}
\numberwithin{defn}{section}

\newcommand{\real}{\mathbb{R}}
\newcommand{\complex}{\mathbb{C}}

\newcommand{\ddd}{\mathbb{D}}

\usepackage{scalerel,stackengine}
\stackMath
\newcommand\reallywidecheck[1]{%
\savestack{\tmpbox}{\stretchto{%
  \scaleto{%
    \scalerel*[\widthof{\ensuremath{#1}}]{\negmedspace\kern0pt\bigwedge\kern0pt}%
    {\rule[-\textheight/2]{1ex}{\textheight}}%WIDTH-LIMITED BIG WEDGE
  }{.5\textheight}% 
}{0.5ex}}%
\stackon[1pt]{#1}{\scalebox{-1}{\tmpbox}}%
}
\newcommand{\rwc}{\reallywidecheck}
\DeclareMathOperator{\supp}{supp}

\DeclareMathOperator{\Tanh}{Tanh}
\DeclareMathOperator{\Sinh}{Sinh}
\DeclareMathOperator{\Cosh}{Cosh}
\DeclareMathOperator{\Sech}{Sech}

\title{Explicit solution of the 1D Schr\"odinger equation}
\author{Peter C.~Gibson}

\date{April 12, 2023} % delete this line to display the current date

\begin{document}

\maketitle

\begin{abstract}
Evaluation of a product integral with values in the Lie group SU(1,1) yields the explicit solution to the impedance form of the Schr\"odinger equation.  Explicit formulas for the transmission coefficient and $S$-matrix of the classical one-dimensional Schr\"odinger operator with arbitrary compactly supported potential are obtained as a consequence.  
The formulas involve operator theoretic analogues of the standard hyperbolic functions, and provide a new window on acoustic and quantum scattering in one dimension. 
\end{abstract}

\begin{center}
MSC 34L40, 34A05, 34L25;
\end{center}

\tableofcontents
%
%\newpage

\section{Introduction}

Inverse theory for the one-dimensional Schr\"odinger operator has been developing since the 1940s. 
Highlights of the vast accumulated literature include the early work of Borg \cite{Bo:1946} and Jost \cite{Jo:1947}, the inverse scattering method of Faddeev \cite{Fa:1964}, the monumental analysis of Deift and Trubowitz \cite{DeTr:1979}, and the inverse spectral theory of Gesztesy and Simon \cite{GeSi:2000}. These all require a certain degree of regularity of the potential function, namely that it be locally integrable.  Moreover, they do not provide a fully explicit form for the solution to the Schr\"odinger equation, or for the reflection coefficient central to scattering theory.  An explicit formula for the reflection coefficient was first established in the preprint \cite{Gi:Ar2021} in the context of a newly developed approach to one-dimensional scattering via singular approximation, which applies to a class of  potential functions of lower regularity than the classical literature. 
The present work extends this approach beyond the reflection coefficient to the full wave field, showing it can be expressed in terms of operator theoretic analogues of the elementary hyperbolic functions.  These operators act on the pre-image of the potential by the Miura map
\[
\alpha\mapsto q=\alpha^2- \alpha^\prime, 
\]
a subtlety which may explain why the formulas presented here were not found earlier. 

The natural algebraic context for our results is the Lie group SU(1,1), consisting of $2\times 2$ complex matrices $U$ of determinant 1 that leave invariant the quadratic form induced by 
\[\Lambda=\begin{pmatrix}1&0\\ 0&-1\end{pmatrix},\] in the sense that  $U^\ast\Lambda U=\Lambda$, or equivalently, $\Lambda U^\ast\Lambda=U^{-1}$.  The corresponding Lie algebra $\mathfrak{su}(1,1)$ consists of $\Lambda$-skew matrices, i.e., matrices $A$ that satisfy $\Lambda A^\ast\Lambda=-A$.  Representation of SU(1,1) is described in Bargmann's foundational paper \cite[\S3b]{Ba:1947}, where it is labelled $\mathfrak{S}_3$; see also \cite[p.201]{Gi:1994}.
This paper's main results may be viewed as an explicit evaluation of a product integral in SU(1,1), with a twist. The approximating finite products are chosen in an unconventional way so as to reveal a key connection to orthogonal polynomials on the unit circle.

In precise technical terms, 
this article concerns bivariate functions $u(x,\sigma)$ governed by the equation 
\begin{equation}\label{wave-fourier}
(\zeta u^\prime)^\prime+\sigma^2\zeta u=0.
\end{equation}
Here $\sigma\in\real$ is frequency, the Fourier dual variable to time; the derivative $u^\prime$ refers to the spatial variable, denoted $x$.  We consider (\ref{wave-fourier}) on a spatial interval $X=(x_0,x_1)$, for some $x_0<x_1$.  The coefficient $\zeta:X\rightarrow\real_+$, called the impedance function, is assumed to be bounded, positive and piecewise absolutely continuous, such that its reciprocal $1/\zeta$ is bounded.  
Equation (\ref{wave-fourier}) is known variously as the Helmholtz equation, or the impedance form of the Schr\"odinger equation.  Its importance stems in part from its well-known connection to the one-dimensional acoustic wave equation 
\begin{equation}\label{wave}
\zeta \partial_{tt}U-\partial_x(\zeta \partial_xU)=0
\end{equation}
and to the classical Schr\"odinger equation
\begin{equation}\label{schrodinger}
-y^{\prime\prime}+qy=\sigma^2y,
\end{equation}
as follows.  Application to (\ref{wave}) of the (distributional) Fourier transform, consistent with the formula
\begin{equation}\label{fourier-transform}
f(x,t)\mapsto \hat{f}(x,\sigma)=\int_{-\infty}^\infty f(x,t)e^{i\sigma t}\,dt,
\end{equation}
yields (\ref{wave-fourier}) with $u=\widehat{U}$.  And if $\zeta^\prime$ is absolutely continuous, then setting 
\begin{equation}\label{schrodinger-cov}
y=\sqrt{\zeta} u\quad\mbox{ and }\quad q=\sqrt{\zeta}^{\,\prime\prime}/\sqrt{\zeta}
\end{equation}
transforms (\ref{wave-fourier}) into (\ref{schrodinger}).  The objective of the present article is to determine the general solution to (\ref{wave-fourier}) in explicit form. In the guise of explicitly formulated (and reasonably simple) operators, this provides a new tool with which to analyze the scattering of acoustic and quantum waves in one dimension. 

We build on \cite{Gi:Ar2021}, which introduces the harmonic exponential operator (see \S\ref{sec-harmonic}), but the present article is self contained.  
The preprint \cite{Gi:Ar2021} is concerned specifically with the reflection coefficient, which is given an explicit formulation.  Here the explicit formulation is extended to the full wave field.  In particular, we obtain a new formula for the transmission coefficient, and for the 
$S$-matrix corresponding to (\ref{schrodinger}). 
The results are significant in that they reduce the analysis of acoustic and quantum scattering to hands-on analysis of a concretely given operators.  High-frequency asymptotics, for example, can simply be read off directly from the given formulas.  

\subsection{Related literature\label{sec-literature}}
The product integral, also known as the multiplicative integral, or ordered matrix exponential, originates in the work of Volterra; see \cite{DoFr:1977}, \cite{Sl:2007} and \cite{GiJo:1990} for background.  The notation for product integration is not standardized; we use that of \cite{Sl:2007}. Previous applications of the product integral to the Schr\"{o}dinger equation have been concerned with asymptotic analysis \cite[\S5.6]{DoFr:1977} and nonlinear equations \cite[Ch.~1, \S2]{FaTa:2007}.  Series formulas in terms of matrices, have been known since the early theory \cite{Pe:1888}.  The present treatment is novel in presenting explicit formulas for individual matrix entries, and relating these to orthogonal polynomials on the unit circle.  See \cite[Ch.~1]{SiOPUC1:2005} and \cite[Ch.~8]{Kh:2008} for background on orthogonal polynomials.  An efficient inverse algorithm for inverse acoustic scattering involving orthogonal polynomials on the unit circle was established in \cite{Gi:JCP2018}.  A distinguishing feature of the present work is the comparatively low regularity of the coefficient $\zeta$ in (\ref{wave-fourier}), corresponding to singularities in the potential $q$ on the level of the derivative of a Dirac delta.  This sets it apart, not only from \cite{GeSi:2000}, but also the perspective of De Branges spaces \cite{BaBePo:2017}, which requires $q\in L^2(a,b)$.  Schr\"{o}dinger potentials of low regularity have been studied in \cite{EcGeNiTe:2013a,EcGeNiTe:2013b}, but not from the perspective of scattering theory; see also \cite{AmRe:2005,FeMeSe:2017}.

\subsection{Organization of the paper\label{sec-organization}}  
Section~\ref{sec-preliminaries} establishes the framework for the paper, addressing: uniqueness of solutions to (\ref{wave-fourier}) (which is proved in Appendix~\ref{sec-appendix-uniqueness}); reformulation of (\ref{wave-fourier}) as a first-order system; the basic jump relation; and definition of the key analytic object, the harmonic exponential operator. Key properties of the harmonic exponential operator are summarized at the end of the section; proofs of these are deferred to \S\ref{sec-analysis}. The main results are presented in Section~\ref{sec-main-results}.  Theorems~\ref{thm-general-solution} and \ref{thm-jumps} in \S\ref{sec-general-solution} formulate the general solution to (\ref{wave-fourier}) in the cases where $\zeta$ is absolutely continuous, and piecewise absolutely continuous, respectively.  Theorem~\ref{thm-schrodinger} in \S\ref{sec-S} describes explicitly the corresponding $S$-matrix for the Schr\"odinger equation.  Section~\ref{sec-applications} applies the main results to analysis of high-frequency asymptotics as well as the structure of the transmission coefficient.  Detailed analysis of the harmonic exponential is carried out in \S\ref{sec-analysis}, providing self-contained proofs of the properties cited earlier in \S\ref{sec-harmonic}.  In particular, the crucial connection to orthogonal polynomials on the unit circle is explained in \S\ref{sec-OPUC}.  Section~\ref{sec-remarks} concludes the paper with a discussion of several aspects of the results.

\section{Preliminaries\label{sec-preliminaries}}

\subsection{Uniqueness of solutions\label{sec-uniqueness}}
Regarding $\sigma$ as a parameter, the precise interpretation of (\ref{wave-fourier}) is as follows.  Solutions to (\ref{wave-fourier}) are functions $u(\cdot,\sigma)\in W^{1,2}(X)$ that satisfy (\ref{wave-fourier}) in the sense of distributions. 
In particular, $u$ and $\zeta u^\prime$ are necessarily continuous in $x$ if $u$ is a solution to (\ref{wave-fourier}).  Solutions in the foregoing sense are uniquely determined, as follows.  
\begin{prop}\label{prop-uniqueness}
Fix $X=(x_0,x_1)$ with $x_0<x_1$, suppose that $\zeta$ is piecewise absolutely continuous on $X$, and let $u(x_0+,\sigma)$ and $u^\prime(x_0+,\sigma)$ be prescribed arbitrarily.  Then for each $\sigma\in\real$, equation (\ref{wave-fourier}) has at most one solution consistent with the given values of $u(x_0+,\sigma)$ and $u^\prime(x_0+,\sigma)$.  
\end{prop}
Proposition~\ref{prop-uniqueness} is proved in Appendix~\ref{sec-appendix-uniqueness}. 

\subsection{Components of the wave field\label{sec-components}}
It is useful to split functions $u(x,\sigma)$ satisfying (\ref{wave-fourier}) into what will be referred to as left-moving and right-moving components (weighted by $\sqrt{\zeta}$), defined as follows.  
For $\sigma\neq 0$, set 
\begin{equation}\label{abA}
\begin{split}
\binom{a}{b}&=\frac{\sqrt{\zeta(x)}}{2}\begin{pmatrix}1&-\frac{1}{i\sigma}\\[5pt] 1&\frac{1}{i\sigma}\end{pmatrix}\binom{u}{u^\prime},\\
A&=\binom{a}{b}.
\end{split}
\end{equation}
Thus, by construction,
\begin{equation}\label{uab}
u=(a+b)/\sqrt{\zeta}\quad\mbox{ and }\quad u^\prime=i\sigma(b-a)/\sqrt{\zeta}.
\end{equation}
If $\zeta$ is constant on some interval, then the inverse Fourier transform of $a(x,\sigma)$ with respect to $\sigma$ is a left-moving wave on that interval. Similarly, $b(x,\sigma)$ corresponds to right-moving waves on intervals where $\zeta$ is constant.  The  terminology of left-moving and right-moving components respectively for $a$ and $b$ is extended to general $\zeta$, although the wave operator no longer factors on intervals where $\zeta$ is non-constant. 

Setting $\alpha=-\frac{1}{2}(\log\zeta)^\prime$ in the case where $\zeta$ is absolutely continuous, equation (\ref{wave-fourier}) is equivalent to the first order system
\begin{equation}\label{wave-system}
A^\prime+\frak{m}A=0\quad\mbox{ where }\quad \frak{m}(x,\sigma)=\begin{pmatrix}i\sigma&\alpha(x)\\ \alpha(x)&-i\sigma\end{pmatrix}. 
\end{equation}
Note that since $\alpha$ is real valued, the coefficient matrix $\frak{m}$ belongs to the Lie algebra $\mathfrak{su}(1,1)$. 

\subsection{The relation between one-sided limits at a point\label{sec-point}}
For $x\in X$, set 
\begin{equation}\label{gamma-J-point}
\gamma(x)=\frac{\zeta(x-)}{\zeta(x+)}\quad\mbox{ and }\quad J=\frac{1}{2\sqrt{\gamma}}\begin{pmatrix}1+\gamma&1-\gamma\\ 1-\gamma&1+\gamma\end{pmatrix}.
\end{equation}
Referring to (\ref{abA}), a simple consequence of continuity of $u$ and $\zeta u^\prime$ at $x$ is the relation
\begin{equation}\label{scattering-at-x}
A(x+,\sigma)=J(x)A(x-,\sigma). 
\end{equation}
Thus $A(\cdot,\sigma)$ is continuous at $x$ if and only $J(x)$ is the identity matrix, which holds if and only if $\zeta$ is continuous at $x$.  Note that $J\in\mbox{SU(1,1)}$.

\subsection{The harmonic exponential operator\label{sec-harmonic}}
The relationship between the coefficient $\alpha$ and the general solution to (\ref{wave-system}) is mediated by what we call the harmonic exponential operator, first introduced in \cite{Gi:Ar2021}. Parameterized by $X=(x_0,x_1)$, the harmonic exponential operator $E^X$ acts on functions $\alpha:Y\rightarrow\real$, where $X\subset Y\subset\real$, such that $\alpha|_X\in L^1(X)$. It is defined by the formula
\begin{equation}\label{hexp-formula}
E^X_\alpha(\sigma):=1+\displaystyle\sum_{j=1}^\infty\thickspace\int\limits_{x_0<s_1<\cdots<s_j<x_1}\negthickspace\negthickspace e^{2i\sigma\kappa(s_1,\ldots,s_j)}\displaystyle\prod_{\nu=1}^j\alpha(s_\nu)\,ds_1\cdots ds_j\qquad(\sigma\in\complex),
\end{equation}
where the term $\kappa(s_1,\ldots,s_j)$ in the exponent depends on the parity of the index $j$:
\begin{multline}\label{kappa}
\kappa(s_1,\ldots,s_j)=\\
\displaystyle\left\{
\begin{array}{cc}
(s_j-s_{j-1})+(s_{j-2}-s_{j-3})+\cdots+(s_2-s_1)&\mbox{ if $j$ is even}\\
(s_j-s_{j-1})+(s_{j-2}-s_{j-3})+\cdots+(s_3-s_2)+(s_1-x_0)&\mbox{ if $j$ is odd}
\end{array}\right..
\end{multline}
The intended interpretation of the notation $E^X_\alpha(\sigma)$ is that $E^X$ is an operator mapping the function $\alpha$ to the complex-valued function $E^X_\alpha$ on $\complex$ (or on $\real$).  But of course $E^X_\alpha(\sigma)$ is also a function of the real parameters $x_0$ and $x_1$, in addition to $\sigma$ and the function $\alpha$.  

In terms of its dependence on $x_1$ with $\sigma$ fixed, $E^X_\alpha(\sigma)$ is absolutely continuous; and direct evaluation from the definition (\ref{hexp-formula}) yields
\begin{equation}\label{x-derivative}
\partial_{x_1}E^{X}_\alpha(\sigma)=e^{2i(x_1-x_0)\sigma}\alpha(x_1)\overline{E^{X}_\alpha(\bar\sigma)}. 
\end{equation}
The function $e^{2i(x_1-x_0)\sigma}\,\overline{E^{X}_\alpha(\bar\sigma)}$ occurring on the right-hand side of (\ref{x-derivative}) is dual to $E^X_\alpha(\sigma)$ in a sense analogous to duality of orthogonal polynomials on the unit circle. Denoting
\begin{equation}\label{dual-E}
\bigl(E_{\alpha}^{X}\bigr)^{\!\ast}(\sigma):=e^{2i(x_1-x_0)\sigma}\,\overline{E^{X}_\alpha(\bar\sigma)}, 
\end{equation}
(\ref{x-derivative}) takes the simple form
\begin{equation}\label{x-derivative-dual}
\partial_{x_1}E^{X}_\alpha=\alpha(x_1)\bigl(E_{\alpha}^{X}\bigr)^{\!\ast}. 
\end{equation}
\begin{remark}\label{remark-duality}
The notion of duality at play in (\ref{dual-E}) applies generally to functions 
\[
F:X\times\complex\rightarrow\complex,
\]
by way of the definition
\[
F^\ast(x,\sigma):=e^{2i(x-x_0)\sigma}\,\overline{F(x,\bar\sigma)}. 
\]
Observe that, according to this definition, duality is involutive: $F^{\ast\ast}=F$.  Based on (\ref{x-derivative-dual}), we call the equation
\begin{equation}\label{harmonic-exponential-equation}
\partial _xF=\alpha F^\ast,\qquad F(x_0+,\sigma)=1
\end{equation}
the \emph{harmonic exponential equation}.  See \S\ref{sec-remarks} for a discussion of its connection to the classical three-term recurrence relation for orthogonal polynomials. 
\end{remark}

Hyperbolic trigonometric operators associated to $E^X$ are defined by analogy with the usual scalar trigonometric functions:
\begin{equation}\label{hyperbolic-operators}
\begin{split}
\Cosh^X_\alpha&=\left(E^X_\alpha+E^X_{-\alpha}\right)/2,\\
\Sinh^X_\alpha&=\left(E^X_\alpha-E^X_{-\alpha}\right)/2,\\
\Sech^X_\alpha&=1/\Cosh^X_\alpha,\\ 
\Tanh^X_\alpha&=\Sinh^X_\alpha/\Cosh^X_\alpha.
\end{split}
\end{equation}

The notation $E^{(x_0,x_1)}$ will be used interchangeably with $E^X$ (and similarly for the other related operators) whenever it is desirable to emphasize the role of the end points of $X$. 
The harmonic exponential operator is useful because it is amenable to direct analysis. We list here some basic properties proved in \S\ref{sec-analysis}. See Propositions~\ref{prop-analytic}\ref{part-entire}, \ref{prop-real-part}, \ref{prop-zero-free}, \ref{prop-lower-bound} and \ref{prop-sech-i} for details. 
\begin{prop}\label{prop-hexp}
The following hold for every $\alpha\in L^1_\real(X)$ and $x\in (x_0,x_1]$. 
\begin{enumerate}[label={(\roman*)},itemindent=0em]
\item With respect to $\sigma\in\complex$, $E^{(x_0,x)}_\alpha(\sigma)$ is entire. \label{part-entire-0}\\[0pt]
\item $\displaystyle\lim_{\stackrel{|\sigma|\rightarrow\infty}{\sigma\in\real}}E^{(x_0,x)}_\alpha(\sigma)=1.$\label{hexp-limit}\\[0pt]
\item For every $\sigma\in\real$, $\Re E^{(x_0,x)}_\alpha(\sigma)\overline{E^{(x_0,x)}_{-\alpha}(\sigma)}=1$.\label{hexp-real-part}\\[0pt]
\item If $\Im\sigma\geq 0$, then
both $E^{(x_0,x)}_\alpha(\sigma)$ and $E^{(x_0,x)}_\alpha(\sigma)+E^{(x_0,x)}_{-\alpha}(\sigma)$ are nonzero. \label{part-zero-free}\\[0pt]
\item If $\Im\sigma\geq0$, then 
$
\bigl|E^{(x_0,x)}_\alpha(\sigma)+E^{(x_0,x)}_{-\alpha}(\sigma)\bigr|\geq\sqrt{2}.
$
\label{part-lower-bound}\\[0pt]
\item $\Sech^{(x_0,x)}_\alpha(i)>0$.\label{part-sech-i}
\end{enumerate}
\end{prop}

\section{Main results\label{sec-main-results}}

\subsection{The general solution\label{sec-general-solution}}

\begin{thm}\label{thm-general-solution}
Suppose $\zeta$ is absolutely continuous on $X=(x_0,x_1)$, and write $\alpha=-\frac{1}{2}(\log\zeta)^\prime$. 
For $x_0<x\leq x_1$ and $\sigma\in\real$, set
\begin{align}\label{zwG}
w(x_0,x,\sigma)&=e^{-i(x-x_0)\sigma}\Sinh_\alpha^{(x_0,x)}(\sigma),\\
z(x_0,x,\sigma)&=e^{-i(x-x_0)\sigma}\Cosh_\alpha^{(x_0,x)}(\sigma),\\
\label{M}
G&=\begin{pmatrix}z&-w\\[5pt] -\bar{w}&\bar{z}\end{pmatrix}.
\end{align}
Then $G(x_0,x,\sigma)\in\mbox{SU(1,1)}$.  
The unique solution to (\ref{wave-system}), parameterized by $A(x_0+,\sigma)$, is 
\begin{equation}\label{general-solution}
A(x,\sigma)=G(x_0,x,\sigma)A(x_0+,\sigma)\qquad(x_0<x<x_1).  
\end{equation}
\end{thm}
\begin{pf} Equation (\ref{x-derivative}) implies
\begin{equation}\label{G-derivative}
\frac{\partial}{\partial x}G(x_0,x,\sigma)=-\begin{pmatrix}i\sigma&\alpha(x)\\ \alpha(x)&-i\sigma\end{pmatrix}G(x_0,x,\sigma).
\end{equation}
Therefore $A(x,\sigma)$ as prescribed by (\ref{general-solution}) satisfies (\ref{wave-system}).  By Proposition~\ref{prop-uniqueness}, the solution is unique.  The definitions of $\Sinh^X$ and $\Cosh^X$ together with Proposition~\ref{prop-hexp}\ref{hexp-real-part} imply 
\[
|z|^2-|w|^2=\Re E^X_\alpha(\sigma)\overline{E^X_{-\alpha}(\sigma)}=1,
\]
proving $G\in\mbox{SU(1,1)}$.  
\end{pf}

\begin{remark}\label{remark-inverse} Set 
\begin{equation}\label{H}
H=G^{-1}=\begin{pmatrix}\bar{z}&w\\ \bar{w}&z\end{pmatrix}.
\end{equation}
Stated in terms of $H$, Theorem~\ref{thm-general-solution} asserts
\begin{equation}\label{H-version}
A(x_0+,\sigma)=H(x_0,x,\sigma)A(x,\sigma)\qquad(x_0<x<x_1). 
\end{equation}
Re-expressing (\ref{G-derivative}) in terms of $H$ yields the equation
\begin{equation}\label{H-derivative}
\frac{\partial}{\partial x}H(x_0,x,\sigma)=H(x_0,x,\sigma)\begin{pmatrix}i\sigma&\alpha(x)\\ \alpha(x)&-i\sigma\end{pmatrix}.
\end{equation}
Section~\ref{sec-product-integral} analyzes the product integral for $H$ rather than $G$. 
\end{remark}

The solution to (\ref{wave-fourier}) on $X$ in terms of $u(x_0+,\sigma)$ and $u^\prime(x_0+,\sigma)$ is 
\begin{equation}\label{u-form}
u(x,\sigma)=u(x_0+,\sigma)\Re\bigl(e^{-i(x-x_0)\sigma}E^{(x_0,x)}_{-\alpha}(\sigma)\bigr)-\frac{1}{\sigma}u^\prime(x_0+,\sigma)\Im\bigl(e^{-i(x-x_0)\sigma}E^{(x_0,x)}_\alpha(\sigma)\bigr)
\end{equation}
Thus if $\binom{u}{u^\prime}$ is real valued anywhere in $X$, it is real valued everywhere. 
Note that the coefficients of $u(x_0+,\sigma)$ and $u^\prime(x_0+,\sigma)$ in (\ref{u-form}) may be alternatively formulated as follows. 
\begin{multline}\label{real-part}
\Re\bigl(e^{-i(x-x_0)\sigma}E^{(x_0,x)}_{-\alpha}(\sigma)\bigr)=\\
\cos(\sigma(x-x_0))+\sum_{j=1}^\infty(-1)^j\rule{-48pt}{0pt}\int\limits_{\stackrel{\rule{0pt}{4pt}}{x_0=s_0<s_1<\cdots<s_j<s_{j+1}=x}}\rule{-48pt}{0pt}\cos\left(\sigma\sum_{\nu=1}^{j+1}(-1)^{j-\nu}(s_\nu-s_{\nu-1})\right)\prod_{\nu=1}^j\alpha(s_j)\,ds_1\cdots ds_j
\end{multline}
\begin{multline}\label{imaginary-part}
\Im\bigl(e^{-i(x-x_0)\sigma}E^{(x_0,x)}_{\alpha}(\sigma)\bigr)=\\
-\sin(\sigma(x-x_0))+\sum_{j=1}^\infty\rule{-42pt}{0pt}\int\limits_{\stackrel{\rule{0pt}{4pt}}{x_0=s_0<s_1<\cdots<s_j<s_{j+1}=x}}\rule{-48pt}{0pt}\sin\left(\sigma\sum_{\nu=1}^{j+1}(-1)^{j-\nu}(s_\nu-s_{\nu-1})\right)\prod_{\nu=1}^j\alpha(s_j)\,ds_1\cdots ds_j
\end{multline} 

Observe by Theorem~\ref{thm-general-solution} that the left-hand limit of the solution at $x_1$ is
\begin{equation}\label{A-minus}
A(x_1-,\sigma)=G(x_1,\sigma)A(x_0+,\sigma),\quad\mbox{ since }\quad G(x_1,\sigma)=G(x_1-,\sigma). 
\end{equation}
Next consider the case where $\zeta$ is piecewise absolutely continuous, with finitely many jump points.  
Theorem~\ref{thm-general-solution} applies to every open interval on which $\zeta$ is absolutely continuous, while the behaviour of $A(x,\sigma)$ across jump points of $\zeta$ is governed by (\ref{scattering-at-x}).  Combining (\ref{scattering-at-x}) with Theorem~\ref{thm-general-solution} therefore yields the following.
\begin{thm}\label{thm-jumps}
Suppose $\zeta$ is piecewise absolutely continuous on $X=(x_0,x_1)$, with jump points confined to $\{y_1,\ldots,y_n\}$, where 
\[
x_0=y_0<y_1<\cdots<y_n<y_{n+1}=x_1. 
\]
Define $J$ as in (\ref{gamma-J-point}), and let $x\in X$ satisfy $y_{j}<x<y_{j+1}$ for some $0\leq j\leq n$.  Then the unique solution to (\ref{wave-fourier}) in terms of the notation (\ref{abA}) is 
\begin{equation}\label{jump-solution}
A(x,\sigma)=G(y_j,x,\sigma)J(y_{j})G(y_{j-1},y_j,\sigma)\cdots J(y_1)G(x_0,y_1,\sigma)A(x_0+,\sigma). 
\end{equation}
\end{thm}

\subsection{The $S$-matrix\label{sec-S}} 
The change of variables (\ref{schrodinger-cov}) yields a potential $q$ two derivatives less regular than $\zeta$.  Therefore, to avoid potentials with singularities on the order of the derivative of a Dirac function, suppose that $\zeta$ is absolutely continuous.  
Suppose further (and without loss of generality) that the support of $\zeta^\prime$ is contained in $X=(x_0,x_1)$.  
Extend $\zeta$ to an absolutely continuous function on $\real$ by setting 
\[
\zeta(x)=\zeta(x_0+)\quad\mbox{ if }\quad x\leq x_0\quad\mbox{ and }\quad\zeta(x)=\zeta(x_1-)\quad\mbox{ if }\quad x\geq x_1. 
\]
The corresponding potential $q$ is then also defined on $\real$, and is zero outside of $X$.   
Theorem~\ref{thm-general-solution} leads easily to an explicit formula for the $S$-matrix corresponding to $q$. 
The classical $S$-matrix is defined in terms of particular solutions $f_1$ and $f_2$ to (\ref{schrodinger}), according to the scheme
\begin{align}
f_1(x,\sigma)&\sim\left\{
\begin{array}{cc}
e^{i\sigma x}&\mbox{ as }x\rightarrow+\infty\\
\frac{1}{T_2(\sigma)}e^{i\sigma x}+\frac{R_2(\sigma)}{T_2(\sigma)}e^{-i\sigma x}&\mbox{ as }x\rightarrow-\infty
\end{array}\right.\label{f1}\\
f_2(x,\sigma)&\sim\left\{
\begin{array}{cc}
\frac{R_1(\sigma)}{T_1(\sigma)}e^{i\sigma x}+\frac{1}{T_1(\sigma)}e^{-i\sigma x}&\mbox{ as }x\rightarrow+\infty\\
e^{-i\sigma x}&\mbox{ as }x\rightarrow-\infty
\end{array}\right.\label{f2}\\
S&=\begin{pmatrix}T_1&R_2\\ R_1&T_2\end{pmatrix}. 
\end{align}
See \cite[\S1]{DeTr:1979}. Note that in the present case, where the support of $q$ is contained in $X$, the asymptotic structure of the functions $f_1,f_2$ for $x\rightarrow+\infty$ is attained as soon as $x\geq x_1$, and for $x\rightarrow-\infty$ when $x\leq x_0$; and $A$ is continuous at both $x_0$ and $x_1$.  Moreover, in the asymptotic regime, the left-moving component has the form $a(x,\sigma)=C_-(\sigma)e^{-i\sigma x}$ and the right-moving component has the form $b(x,\sigma)=C_+(\sigma)e^{i\sigma x}$, for appropriate coefficients $C_\pm$.  Explicit formulas for the entries of the $S$-matrix obtain from these facts, as follows. 

Consider first $f_1$.  In this case, according to (\ref{f1}), 
\begin{equation}\label{A0A1}
A(x_0,\sigma)=\binom{\frac{R_2}{T_2}e^{-i\sigma x_0}}{\frac{1}{T_2}e^{i\sigma x_0}}\quad\mbox{ and }\quad A(x_1,\sigma)=\binom{0}{e^{i\sigma x_1}}.
\end{equation}
Given continuity of $A$ at $x_0$ and $x_1$, Theorem~\ref{thm-general-solution} asserts 
\begin{equation}\label{application-of-thm1}
A(x_1,\sigma)=G(x_0,x_1,\sigma)A(x_0,\sigma),
\end{equation}
allowing one to solve for $R_2$ and $T_2$ in terms of $G(x_0,x_1,\sigma)$ using (\ref{A0A1}).  The result is
\begin{equation}\label{R2T2}
R_2(\sigma)=e^{2i\sigma x_0}\Tanh_\alpha^X(\sigma)\quad\mbox{ and }\quad T_2(\sigma)=\Sech^X_\alpha(\sigma).  
\end{equation}
A similar analysis applied to $f_2$ yields
\begin{equation}\label{R1T1}
R_1(\sigma)
=-e^{-2i\sigma x_0}\frac{\overline{\Sinh^X_\alpha(\sigma)}}{\Cosh^X_\alpha(\sigma)}\quad\mbox{ and }\quad T_1(\sigma)=\Sech^X_\alpha(\sigma).  
\end{equation}
In summary, we have proved the following. 
\begin{thm}\label{thm-schrodinger}
Let $q=\bigl(\sqrt{\zeta}\bigr)^{\prime\prime}/\sqrt{\zeta}$, where $\zeta:\real\rightarrow\real_{>0}$ is such that $\supp\zeta^\prime\subset X=(x_0,x_1)$.  Then the $S$-matrix for the classical Schrodinger equation (\ref{schrodinger}) is given by the formula
\begin{equation}\label{S}
S=\Sech^X_\alpha\begin{pmatrix}1&e^{2i\sigma x_0}\Sinh^X_\alpha\\[5pt] -e^{-2i\sigma x_0}\overline{\Sinh^X_\alpha}&1\end{pmatrix}
\end{equation}
where $\alpha=-\frac{1}{2}(\log\zeta)^\prime$.
\end{thm}
The $S$-matrix (\ref{S}) is easily seen to be unitary directly from the given formulas (in keeping with established theory).  What is new here is the explicit form in terms of the harmonic exponential operator, and in particular the representation of the transmission coefficient as the hyperbolic secant operator.  

Note that $q=\bigl(\sqrt{\zeta}\bigr)^{\prime\prime}/\sqrt{\zeta}$ and $\alpha=-\frac{1}{2}(\log\zeta)^\prime$ are related by the Miura map,
\begin{equation}\label{miura}
q=\alpha^2-\alpha^\prime. 
\end{equation}
Kappeler et al.~have shown that a potential has the form (\ref{miura}) 
if and only if the corresponding Schrodinger operator has no bound states \cite[Thm.~1.1]{KaPeShTo:2005}.  Thus Theorem~\ref{thm-schrodinger} applies only to potentials without bound states; but this is not really a restriction, as Deift and Trubowitz explain in \cite[\S3]{DeTr:1979}. 
Given a compactly supported potential $\tilde{q}$ with bound states, there is an explicitly determined associated potential $q$ obtained by removing the bound states, with a corresponding $\alpha$ related to $q$ by (\ref{miura}). Multiplication of the entries of (\ref{S}) by appropriate Blaschke products then gives the $S$-matrix for $\tilde{q}$ in explicit form.  Thus the results of Deift and Trubowitz and of Kappeler et al.~effectively extend Theorem~\ref{thm-schrodinger} to arbitrary compactly supported potential functions.  See \S\ref{sec-remarks} for further discussion of the Miura map.

\section{Applications\label{sec-applications}}
\subsection{High-frequency asymptotics\label{sec-high-frequency}}
Proposition~\ref{prop-hexp}\ref{hexp-limit} allows one to read off the high-frequency asymptotics of $A$ directly from Theorem~\ref{thm-jumps} in the case where $\zeta$ is piecewise absolutely continuous. 
 
First consider the basic case of a step function $\zeta_0$ on $X=(x_0,x_1)$, with jump points confined to $\{y_1,\ldots,y_n\}$, where 
\[
x_0=y_0<y_1<\cdots<y_n<y_{n+1}=x_1. 
\]
Since $\zeta_0$ is constant on intervals $(y_{j-1},y_j)$, the solution to (\ref{wave-system}) with $\zeta=\zeta_0$ is
\[
A(x,\sigma)=\begin{pmatrix}e^{-i(x-y_{j-1})\sigma}&0\\ 0&e^{i(x-y_{j-1})\sigma}\end{pmatrix}A(y_{j-1}+,\sigma),
\]
for $x\in(y_{j-1},y_j)$, for each $1\leq j\leq n+1$.  By (\ref{scattering-at-x}), it follows that $A(y_j+,\sigma)$ relates to $A(y_{j-1}+,\sigma)$ according to the formula
\[
A(y_j+,\sigma)=J(y_j)D_j(\sigma)A(y_{j-1}+,\sigma)\quad\mbox{ where }\quad D_j(\sigma)=\begin{pmatrix}e^{-i(y_j-y_{j-1})\sigma}&0\\ 0&e^{i(y_j-y_{j-1})\sigma}\end{pmatrix}.
\]
Any $x\in X$ which is not a jump point of $\zeta_0$ satisfies $y_{j}<x<y_{j+1}$ for some $0\leq j\leq n$. The above formula easily yields the solution to (\ref{wave-system}) at $x$ in terms of $A(x_0+,\sigma)$,
\begin{equation}\label{step-solution}
A(x,\sigma)=\begin{pmatrix}e^{-i(x-y_j)\sigma}&0\\ 0&e^{i(x-y_j)\sigma}\end{pmatrix}J(y_j)D_j(\sigma)\cdots J(y_1)D_1(\sigma)A(x_0+,\sigma).
\end{equation}

Now consider the case where $\zeta$ is piecewise absolutely continuous, as in the hypothesis of Theorem~\ref{thm-jumps}.  Note that a piecewise absolutely continuous step function $\zeta$ on $X$ has a unique factorization of the form
\begin{equation}\label{factorization}
\zeta=\zeta_{0\,}\zeta_1
\end{equation}
where $\zeta_0$ is a step function with $\zeta_0(x_0+)=1$, and $\zeta_1$ is absolutely continuous.
Applying Proposition~\ref{prop-hexp}\ref{hexp-limit} to the formulas (\ref{hyperbolic-operators}) yields
\[
\Cosh^X_\alpha(\sigma)\sim 1\quad\mbox{ and }\quad \Sinh^X_\alpha(\sigma)\sim0\quad\mbox{ as }|\sigma|\rightarrow\infty. 
\]
Therefore the matrices $G(y_{j-1},y_j,\sigma)$ occurring in Theorem~\ref{thm-jumps} have the asymptotic structure
\begin{equation}\label{G-asymptotic}
G(y_{j-1},y_j,\sigma)\sim\begin{pmatrix}e^{-i(y_j-y_{j-1})\sigma}&0\\ 0&e^{i(y_j-y_{j-1})\sigma}\end{pmatrix}=D_j(\sigma)\quad\mbox{ as }\quad |\sigma|\rightarrow\infty. 
\end{equation}
The matrices $J(y_j)$ are of course independent of $\sigma$.  Replacing each $G(y_{j-1},y_j,\sigma)$ in (\ref{jump-solution}) by its high-frequency asymptote (\ref{G-asymptotic}) gives the asymptotic formula for $A(x,\sigma)$.  
The resulting formula is precisely (\ref{step-solution}).  This shows that the high-frequency asymptotics of the solution to (\ref{wave-fourier}) depend solely on the step function factor $\zeta_0$ in (\ref{factorization}); they are completely independent of the absolutely continuous factor $\zeta_1$.  Moreover, since the jump matrices $J(y_j)$ are determined by the ratios $\zeta(y_j-)/\zeta(y_j+)$, it is the relative change in impedance at the jumps, rather than the absolute change, that determines the high-frequency asymptotics of the solution to the wave equation.

\subsection{The transmission coefficient is an outer function\label{sec-transmission}}
Suppose a layered medium is represented by an absolutely continuous impedance function $\zeta:\real\rightarrow\real_+$ and governed by the wave equation (\ref{wave}), where $\zeta$ is constant outside of $X=(x_0,x_1)$.  A currently open problem in acoustics is to constrain the unknown structure of $\zeta$ on $X$ by measuring the transmission of sound across $X$. An equivalent problem is to determine all impedance functions consistent with phaseless reflection data; see \cite{Ch:2022}.
In the ideal transmission experiment, a Dirac pulse is transmitted toward $x_0$ from the left, and the resulting wave field $\tau(t):=U(x_1,t)$ is measured at $x_1$, for $t>0$.  But in practice, instrumentation is often designed to measure instead the power spectrum of $\tau$, by transmitting a unit-amplitude sinusoidal signal toward $x_0$ and measuring the steady state amplitude of the resulting transmitted signal at $x_1$, over a range of frequencies. It follows from Theorem~\ref{thm-general-solution} that the Fourier transform of idealized transmission data is
\[
\hat\tau(\sigma)=e^{i(x_1-x_0)\sigma}\Sech^X_\alpha(\sigma),
\]
where $\alpha=-\frac{1}{2}(\log\zeta)^\prime$.  This idealized transmission data can be related to measurement of the power spectrum using Proposition~\ref{prop-hexp}, as follows. 
\begin{thm}\label{thm-sech-outer}
For every $\alpha\in L^1_\real(X)$, the function $\Sech^X_\alpha(\sigma)$ is bounded and holomorphic on the upper half plane $\Im\sigma\geq0$, and hence outer.  The formula
\begin{equation}\label{sech-representation}
\Sech^X_\alpha(\sigma)=\exp\left(\frac{1}{i\pi}\int_{\real}\frac{1+s\sigma}{s-\sigma}\log\bigl|\Sech^X_\alpha(s)\bigr|\,\frac{ds}{1+s^2}\right)\qquad(\Im\sigma\geq0)
\end{equation}
expresses $\Sech^X_\alpha$ in terms of the restriction of its modulus to $\real$. 
\end{thm}
\begin{pf}
Parts~\ref{part-entire} and \ref{part-zero-free} of Proposition~\ref{prop-hexp} imply $\Sech^X_\alpha$ is holomorphic on the upper half plane. It is bounded by Proposition~\ref{prop-hexp}\ref{part-lower-bound}.  Thus $\Sech^X_\alpha$ is outer.  The standard representation for outer functions then yields (\ref{sech-representation}), provided $\Sech^X_\alpha(i)>0$.  See \cite[p.~133]{Ho:1988}. Proposition~\ref{prop-hexp}\ref{part-sech-i} guarantees $\Sech^X_\alpha(i)>0$.\end{pf}

Assuming the thickness $x_1-x_0$ of the interval $X$ is known, Theorem~\ref{thm-sech-outer} implies \emph{the transmission data $\tau$ is completely determined by its power spectrum},
\[
|\hat\tau(\sigma)|^2=\bigl|\Sech^X_\alpha(\sigma)\bigr|^2. 
\]
So, in principle, the technically more difficult, and less robust, time-domain experiment required to measure $\tau$ directly does not provide any additional information about $\zeta$ beyond what is already implicit in the power spectrum.

\section{Analysis of the harmonic exponential\label{sec-analysis}}
Here $L^1_\real(X)$ denotes real-valued Lebesgue integrable functions on $X$.  Let $\mathscr{A}$ denote the algebra of entire functions (i.e.~holomorphic functions $f:\complex\rightarrow\complex$), endowed with the standard topology of uniform convergence on compact sets.  

\subsection{Properties of the harmonic exponential\label{sec-properties-1}}
For convenience we repeat here the formula (\ref{hexp-formula}),
\[
E^X_\alpha(\sigma):=1+\displaystyle\sum_{j=1}^\infty\thickspace\int\limits_{x_0<s_1<\cdots<s_j<x_1}\negthickspace\negthickspace e^{2i\sigma\kappa(s_1,\ldots,s_j)}\displaystyle\prod_{\nu=1}^j\alpha(s_\nu)\,ds_1\cdots ds_j\qquad(\sigma\in\complex),
\]
with some additional notation, as follows. For each $j\geq1$ and permutation $\pi$ of the index set $\{1,\ldots,j\}$, denote by $\mathscr{S}^X_{j,\pi}$ the simplex 
\begin{equation}\label{simplex}
\mathscr{S}^X_{j,\pi}=\left\{(s_1,\ldots,s_j)\in\real^j\,\left|\,x_0<s_{\pi(1)}<\cdots<s_{\pi(j)}<x_1\right.\right\}. 
\end{equation}
The $j!$ simplices (\ref{simplex}) are pairwise disjoint, congruent to one another, and cover the  $j$-cube $(x_0,x_1)^j$ up to a set of measure zero. Denote by $C_j$ the $j$th summand in the series defining $E^X_\alpha$,
\begin{equation}\label{Cj}
C_j(\sigma)=\int\limits_{x_0<s_1<\cdots<s_j<x_1}\negthickspace\negthickspace e^{2i\sigma\kappa(s_1,\ldots,s_j)}\displaystyle\prod_{\nu=1}^j\alpha(s_\nu)\,ds_1\cdots ds_j. 
\end{equation}
Observe by (\ref{kappa}) that for $\sigma=s+it$ $(s,t\in\real)$, 
\begin{equation}\label{beta-sigma}
\left|e^{2i\sigma\kappa(s_1,\ldots,s_j)}\right|\leq\beta_\sigma:=\left\{
\begin{array}{cc}
1&\mbox{ if } t\geq 0\\
e^{2|t|(x_1-x_0)}&\mbox{ if } t<0
\end{array}\right..
\end{equation}
Thus, 
\begin{equation}\label{Cj-bound}
\begin{split}
\left|C_j(\sigma)\right|&\leq\beta_\sigma\int\limits_{\mathscr{S}^X_{j,\rm{id}}}\prod_{\nu=1}^j|\alpha(s_\nu)|\,ds_1\cdots ds_j\\
&=\frac{\beta_\sigma}{j!}\sum_{\pi}\int\limits_{\mathscr{S}^X_{j,\pi}}\prod_{\nu=1}^j|\alpha(s_\nu)|\,ds_1\cdots ds_j\\
&=\frac{\beta_\sigma}{j!}\int\limits_{(x_0,x_1)^j}\prod_{\nu=1}^j|\alpha(s_\nu)|\,ds_1\cdots ds_j\\
&=\frac{\beta_\sigma}{j!}\left(\int_{x_0}^{x_1}|\alpha(s)|\,ds\right)^j\\
&=\beta_\sigma\|\alpha\|^j/j!.\\
\end{split}
\end{equation}
The estimate (\ref{Cj-bound}) leads easily to the following. 
\begin{prop}\label{prop-analytic}Values of the harmonic exponential operator are entire, and bounded in the upper half plane. More precisely:
\begin{enumerate}[label={(\roman*)},itemindent=0em]
\item\label{part-entire}
$E^X_\alpha\in\mathscr{A}$ for every $\alpha\in L^1_\real(X)$;
\item\label{part-bounded}
if $\alpha\in L^1_\real(X)$ and $\Im\sigma\geq 1$, then $\bigl|E^X_\alpha(\sigma)\bigr|\leq e^{\|\alpha\|}$.
\end{enumerate}
\end{prop}
\begin{pf}
As a function of $\sigma$, each $C_j(\sigma)$ as defined in (\ref{Cj}) is entire. It follows from (\ref{Cj}) that the series defining $E^X_\alpha(\sigma)$ converges absolutely, uniformly on compact sets.  Therefore $E^X_\alpha$ is entire, proving \ref{part-entire}. Part~\ref{part-bounded} follows immediately from (\ref{Cj-bound}). 
\end{pf}
\begin{prop}\label{prop-limit-1}
For any $\alpha\in L^1_\real(X)$, 
$
\displaystyle\lim_{\stackrel{|\sigma|\rightarrow\infty}{\sigma\in\real}}E^X_\alpha(\sigma)=1. 
$
\end{prop}
\begin{pf}
Set $h=E^X_\alpha-1$. Fix $j\geq 1$, write $s=(s_1,\ldots,s_j)$, and note by (\ref{kappa}) that $0\leq\kappa(s)\leq (x_1-x_0)$ for $s\in\mathscr{S}^X_{j,\rm{id}}$. For $0\leq t\leq 2(x_1-x_0)$, denote by $\mathcal{H}_t$ the hyperplane
\[
\mathcal{H}_t=\left\{\left.s\in\real^j\,\right|\,2\kappa(s)=t\right\},
\]
and let $\mu(s)$ denote ordinary Lebesgue measure on $\mathcal{H}_t$.  Re-write $C_j(\sigma)$ as 
\[
\begin{split}
C_j(\sigma)&=\frac{1}{2\sqrt{j}}\int_{t=0}^{2(x_1-x_0)}\int_{\mathscr{S}^X_{j,\rm{id}}\cap\mathcal{H}_t}e^{2i\sigma\kappa(s_1,\ldots,s_j)}\prod_{\nu=1}^j\alpha(s_\nu)\,d\mu(s)dt \\
&=\int_{t=-\infty}^{\infty}\left(\frac{1}{2\sqrt{j}}\chi_{(0,2(x_1-x_0))}(t)\int_{\mathscr{S}^X_{j,\rm{id}}\cap\mathcal{H}_t}\prod_{\nu=1}^j\alpha(s_\nu)\,d\mu(s)\right)e^{i\sigma t}dt \\
\end{split}
\]
and take the inverse Fourier transform with respect to $\sigma\in\real$ to yield
\begin{equation}\label{A-j-Fourier}
\rwc{C_j}(t)=\frac{1}{2\sqrt{j}}\chi_{(0,2(x_1-x_0))}(t)\int_{\mathscr{S}^X_{j,\rm{id}}\cap\mathcal{H}_t}\prod_{\nu=1}^j\alpha(s_\nu)\,d\mu(s),
\end{equation}
the interpretation of which in the case $j=1$ is 
\begin{equation}\label{A-Fourier}
\rwc{C_j}(t)=\frac{1}{2}\chi_{(0,2(x_1-x_0))}(t)\alpha\bigl(\textstyle\frac{1}{2}t+x_0\bigr).
\end{equation}
The formulations (\ref{A-j-Fourier}) and (\ref{A-Fourier}) make clear that $\supp\rwc{C_j}\subset[0,2(x_1-x_0)]$
for each $j\geq 1$. 
By definition,
\[
h(\sigma)=E^{(x_0,y)}_\alpha(\sigma)-1=\sum_{j\geq 1}C_j(\sigma).
\]
Therefore 
\[
\supp\rwc{h}\subset\bigcup_{j\geq 1}\supp\rwc{C_j}\subset[0,2(x_1-x_0)].
\]
Furthermore,
\begin{equation*}
\begin{split}
\int_{\real}\left|\rwc{C_j}(t)\right|dt&\leq \frac{1}{2\sqrt{j}}\int\limits_0^{2(x_1-x_0)}\int\limits_{\mathscr{S}^X_{j,\rm{id}}\cap\mathcal{H}_t}\prod_{\nu=1}^j\left|\alpha(s_\nu)\right|\,d\mu(s)\,dt\\
&=\int\limits_{\mathscr{S}^X_{j,\rm{id}}}\prod_{\nu=1}^j\left|\alpha(s_\nu)\right|\,ds\\
&=\frac{1}{j!}\int\limits_{(x_0,x_1)^j}\prod_{\nu=1}^j\left|\alpha(s_\nu)\right|\,ds\\
&=\frac{1}{j!}\left(\int_{x_0}^{x_1}\left|\alpha\right|\right)^j\\
&= \|\alpha\|^j/j!.
\end{split}
\end{equation*}
It follows by Minkowski's inequality that 
\begin{equation*}
\int_\real\bigl|\rwc{h}\bigr|\leq\sum_{j\geq 1}\int_{\real}\left|\rwc{C_j}\right|\leq\sum_{j\geq1}\|\alpha\|^j/j!=e^{\|\alpha\|}-1,
\end{equation*}
proving $\rwc{h}\in L^1(\real)$.  The Riemann-Lebesgue lemma therefore implies $h(\sigma)\rightarrow0$ as $|\sigma|\rightarrow\infty$ $(\sigma\in\real)$.  The desired result follows. 
\end{pf}
\begin{lem}\label{lem-inequality}
Let $\alpha,\tilde\alpha\in L^1_\real(X)$. Then
\[
\left|E^X_{\tilde\alpha}(\sigma)-E^X_\alpha(\sigma)\right|\leq\beta_\sigma\|\tilde\alpha-\alpha\|e^{\max\{\|\alpha\|,\|\tilde\alpha\|\}}.
\]
\end{lem}
\begin{pf}
Set $\varepsilon=\tilde\alpha-\alpha$, and let $C_j$ and $\widetilde{C}_j$ denote the respective $j$th summands of $E^X_\alpha$ and $E^X_{\tilde\alpha}$. Then
\begin{equation}\nonumber
\begin{split}
\left|\widetilde{C}_j(\sigma)-C_j(\sigma)\right|&=\left|\int\limits_{\mathscr{S}^X_{j,\rm{id}}}e^{2i\kappa(s_1,\ldots,s_j)}\left(\prod_{\nu=1}^j\tilde\alpha(s_\nu)-\prod_{\nu=1}^j\alpha(s_\nu)\right)\,ds_1\cdots ds_j\right|\\
&\leq\beta_\sigma\int\limits_{\mathscr{S}^X_{j,\rm{id}}}\left|\prod_{\nu=1}^j\tilde\alpha(s_\nu)-\prod_{\nu=1}^j\alpha(s_\nu)\right|\,ds_1\cdots ds_j\\
&=\frac{\beta_\sigma}{j!}\int\limits_{(x_0,x_1)^j}\left|\prod_{\nu=1}^j\bigl(\alpha(s_\nu)+\varepsilon(s_\nu)\bigr)-\prod_{\nu=1}^j\alpha(s_\nu)\right|\,ds_1\cdots ds_j\\
&\leq\frac{\beta_\sigma}{(j-1)!}\max\left\{\int_{x_0}^{x_1}|\alpha+\varepsilon|,\int_{x_0}^{x_1}|\alpha|\right\}^{j-1}\int_{x_0}^{x_1}|\varepsilon|\\
&=\beta_\sigma\|\tilde\alpha-\alpha\|\frac{\max\{\|\alpha\|,\|\tilde\alpha\|\}^{j-1}}{(j-1)!}.
\end{split}
\end{equation}
Since
\[
\left|E^X_{\tilde\alpha}(\sigma)-E^X_\alpha(\sigma)\right|\leq\sum_{j=1}^\infty\left|\widetilde{C}_j(\sigma)-C_j(\sigma)\right|,
\]
the result follows. 
\end{pf}

It follows from Proposition~\ref{prop-limit-1} that $E^X_\alpha|_\real\in L^\infty(\real)$. To reduce notational clutter, the same symbol will be used for $E^X_\alpha$ and its restriction to the real line; the intended interpretation should be clear from context.  
\begin{prop}\label{prop-continuity}
The harmonic exponential operator is continuous, whether one considers its values as functions on $\complex$ or $\real$, as follows. 
\begin{enumerate}[label={(\roman*)},itemindent=0em]
\item\label{part-complex-continuity}
The mapping $E^X:L^1_\real(X)\rightarrow\mathscr{A}$ is continuous. 
\item\label{part-real-continuity}
The mapping $E^X:L^1_\real(X)\rightarrow L^\infty(\real)$ is continuous.
\end{enumerate}
\end{prop}
\begin{pf}
Suppose $\alpha_n\rightarrow\alpha$ in $L^1_\real(X)$.  By Lemma~\ref{lem-inequality}, 
\[
\left|E^X_{\alpha_n}(\sigma)-E^X_\alpha(\sigma)\right|\leq\beta_\sigma\|\alpha_n-\alpha\|e^{\max\{\|\alpha_n\|,\|\tilde\alpha\|\}},
\]
so $E^X_{\alpha_n}\rightarrow E^X_\alpha$ uniformly on compact subsets of $\complex$, proving \ref{part-complex-continuity}. Since $\beta_\sigma=1$ if $\sigma\in\real$, convergence is uniform also on $\real$, proving \ref{part-real-continuity}. 
\end{pf}

\subsection{The matrix $H$ as a product integral\label{sec-product-integral}}
Referring to Theorem~\ref{thm-general-solution}, we relate the product integral formulation of 
\[
H=G^{-1}=\begin{pmatrix}\bar{z}&w\\ \bar{w}&z\end{pmatrix}
\] 
to orthogonal polynomials on the unit circle.  A first step is to establish some notation. Fix $\alpha\in L^1_{\real}(X)$ and $x\in(x_0,x_1]$. For each $n\geq 1$, set 
\begin{equation}\label{partitions}
\Delta_n=(x-x_0)/n\quad\mbox{ and }\quad y_{n,j}=x_0+j\Delta_n\quad(0\leq j\leq n),
\end{equation}
so that 
\[
x_0=y_{n,0}<y_{n,1}<\cdots<y_{n,n-1}<y_{n,n}=x
\]
is an equally-spaced partition of the interval $(x_0,x)$.  
\begin{remark}\label{remark-x-dependency}
The points $y_{n,j}$ defined in (\ref{partitions}) depend on $x$.  Although not explicitly indicated in the notation, this dependence on $x$ is to be understood throughout the present section, whenever the notation $y_{n,j}$ is invoked. 
\end{remark}
Denote by $M_0$ the matrix
\begin{equation}\label{M0}
M_0=\begin{pmatrix}1&1\\ -1&1\end{pmatrix}.
\end{equation}
For $\sigma\in\real$, observe that 
\begin{equation}\label{observe-H}
e^{i(x-x_0)\sigma}M_0H(x_0,x,\sigma)=
\begin{pmatrix}
e^{2i(x-x_0)\sigma}\overline{E^{(x_0,x)}_\alpha(\bar\sigma)}&E^{(x_0,x)}_\alpha(\sigma)\\ 
-e^{2i(x-x_0)\sigma}\overline{E^{(x_0,x)}_{-\alpha}(\bar\sigma)}&E^{(x_0,x)}_{-\alpha}(\sigma)
\end{pmatrix}.
\end{equation}
As in (\ref{wave-system}), denote
\[
\frak{m}(y,\sigma)=\begin{pmatrix}i\sigma&\alpha(y)\\ \alpha(y)&-i\sigma\end{pmatrix}.
\]

If $\alpha$ is continuous, by standard results concerning product integration (see \cite[Thm.~2.5]{Sl:2007} and \cite[Thm.~2.1]{DoFr:1977}), equation (\ref{H-derivative}) implies 
\begin{multline}\label{H-limit-1}
H(x_0,x,\sigma)=\\ \left(I+\mathfrak{m}(y,\sigma)\,dy\right)\prod_{x_0}^x=\lim_{n\rightarrow\infty}\left(I+\Delta_n\mathfrak{m}(y_{n,1},\sigma)\right)\cdots\left(I+\Delta_n\mathfrak{m}(y_{n,n},\sigma)\right).
\end{multline}
This is a standard form of the product integral, expressed as a limit of ``Riemann products."  A key insight of the present paper is to notice that a slight  perturbation the factors in the standard Riemann product reveals a useful connection to orthogonal polynomials.  We first sketch the appropriate perturbation, then rigorously justify the inherent approximations afterward.

The factors $I+\Delta_n\mathfrak{m}(y_{n,j},\sigma)$ have the form
\begin{equation}\label{matrix-perturbation}
I+\Delta_n\mathfrak{m}(y_{n,j},\sigma)=
\begin{pmatrix}1+i\Delta_n\sigma&\Delta_n\alpha(y_{n,j})\\
\Delta_n\alpha(y_{n,j})&1-i\Delta_n\sigma\end{pmatrix}\cong
\begin{pmatrix}e^{i\Delta_n\sigma}&\Delta_n\alpha(y_{n,j})\\
\Delta_n\alpha(y_{n,j})&\rule{0pt}{13pt}e^{-i\Delta_n\sigma}\end{pmatrix},
\end{equation}
the approximation being accurate to $O(\Delta_n^2)$. The right-hand matrix in (\ref{matrix-perturbation}) has two advantages over that on the left: it belongs to $SU(1,1)$; and products of such matrices yield orthogonal polynomials on the unit circle, as follows. Set
\begin{equation}\label{Mnj}
z=e^{2i\Delta_n\sigma},\quad r_{n,j}=\Delta_n\alpha(y_{n,j}),\quad
M_{n,j}=\begin{pmatrix}z&z\,r_{n,j}\\ r_{n,j}&1\end{pmatrix},\quad 
D_n=\begin{pmatrix}1&0\\ 0&e^{i\Delta_n\sigma}\end{pmatrix}, 
\end{equation}
so that 
\begin{equation}\label{Mnj-property}
\begin{pmatrix}e^{i\Delta_n\sigma}&\Delta_n\alpha(y_{n,j})\\
\Delta_n\alpha(y_{n,j})&\rule{0pt}{13pt}e^{-i\Delta_n\sigma}\end{pmatrix}
=e^{-i\Delta_n\sigma}D_nM_{n,j}D_n^{-1}.
\end{equation}
Using the approximation (\ref{matrix-perturbation}), the product-integral formulation of $H$ (\ref{H-limit-1}) yields 
\[ 
H(x_0,x,\sigma)\cong e^{-i(x-x_0)\sigma}D_nM_{n,1}\cdots M_{n,n}D_n^{-1}.
\] 
If $\Delta_n$ is small relative to $\sigma$, then $D_n\cong I$, and the latter approximation simplifies to
\begin{equation}\label{H-approximation-2}
H(x_0,x,\sigma)\cong e^{-i(x-x_0)\sigma}M_{n,1}\cdots M_{n,n},
\end{equation}
which in turn by (\ref{observe-H}) yields 
\begin{equation}\label{OPUC-1}
M_0M_{n,1}\cdots M_{n,n}\cong\begin{pmatrix}
e^{2i(x-x_0)\sigma}\overline{E^{(x_0,x)}_\alpha(\bar\sigma)}&E^{(x_0,x)}_\alpha(\sigma)\\ 
-e^{2i(x-x_0)\sigma}\overline{E^{(x_0,x)}_{-\alpha}(\bar\sigma)}&E^{(x_0,x)}_{-\alpha}(\sigma)
\end{pmatrix}.
\end{equation}
The significance of (\ref{OPUC-1}) is that the product $M_0M_{n,1}\cdots M_{n,n}$ is a classical formula for orthogonal polynomials on the unit circle.  

\subsection{Orthogonal polynomials on the unit circle\label{sec-OPUC}} We first summarize some well-known facts concerning orthogonal polynomials on the unit circle.  See \cite[Ch.~1]{SiOPUC1:2005} and \cite[Ch.~8]{Kh:2008}.  Keeping the notation of the previous section, set 
\begin{equation}\label{OPUC-notation}
\begin{pmatrix}
\Psi_{n,j}(z)&\Psi^\ast_{n,j}(z)\\ -\Phi_{n,j}(z)&\Phi_{n,j}^\ast(z)
\end{pmatrix}
:=M_0M_{n,1}\cdots M_{n,j}\qquad(0\leq j\leq n). 
\end{equation}
Suppose $|r_{n,j}|<1$ $(1\leq j\leq n)$; if $\alpha$ is bounded, this condition is fulfilled automatically once $n$ is sufficiently large. Then $\Phi_{n,j}(z)$ is a monic polynomial in $z$ of degree $j$ $(0\leq j\leq n)$, and the sequence $\Phi_{n,0},\ldots,\Phi_{n,n}$ is orthogonal on the unit circle with respect to the probability measure 
\begin{equation}\label{mun}
d\mu_n\bigl(e^{i\theta}\bigr)=\frac{\prod_{j=1}^n(1-| r_{n,j}|^2)}{2\pi\left|\Phi_{n,n}(e^{i\theta})\right|^2}\,d\theta\qquad (-\pi<\theta\leq \pi).
\end{equation}
The polynomial $\Phi_{n,j}$ and its dual $\Phi_{n,j}^\ast$ are related by the equation
\begin{equation}\label{dual-phi}
\Phi^\ast_{n,j}(z)=z^j\overline{\Phi_{n,j}(1/\bar{z})}\quad\mbox{ or equivalently, }\quad\Phi_{n,j}=z^j\overline{\Phi_{n,j}^\ast(1/\bar{z})}\qquad(0\leq j\leq n). 
\end{equation}
Thus $\Phi_{n,j}$ is a polynomial of degree at most $j$, but not monic. The classical recurrence
\begin{equation}\label{recurrence}
\Phi_{n,j+1}(z)=z\Phi_{n,j}(z)- r_{n,j+1}\Phi_{n,j}^\ast(z)\qquad(0\leq j\leq n-1)
\end{equation}
is an immediate consequence of the matrix formulation (\ref{OPUC-notation}), as is the recurrence
\begin{equation}\label{recurrence-psi-ast}
\Psi^\ast_{n,j+1}(z)=\Psi_{n,j}^\ast(z)+zr_{n,j+1}\Psi_{n,j}(z)\qquad(0\leq j\leq n-1).
\end{equation} 
In general, all the same formulas hold with $\Psi_{n,j}$ in place of $\Phi_{n,j}$, provided each occurrence of $ r_{n,j}$ is replaced by $- r_{n,j}$. Indeed, replacing $\alpha$ by $-\alpha$ (and hence replacing each $ r_{n,j}$ by its negative) corresponds to interchanging each $\Psi_{n,j}$ with $\Phi_{n,j}$, and $\Psi_{n,j}^\ast$ with $\Phi_{n,j}^\ast$. 
For later reference, we note one further known fact concerning orthogonal polynomials on the unit circle. 
\begin{prop}[{\cite[Thm.~1.7.2]{SiOPUC1:2005},\cite[Prop.~4.2]{Gi:JAT2023}}]\label{prop-polynomials-zero-free}
The polynomials 
\[\Phi_{n,j}^\ast,\quad \Psi_{n,j}^\ast\quad\mbox{ and }\quad\Phi_{n,j}^\ast+\Psi_{n,j}^\ast\] 
are zero free on the closed unit disk $\overline{\ddd}$ $(0\leq j\leq n)$. 
\end{prop}
The entries of the matrix (\ref{OPUC-notation}) can be worked out explicitly. Observe that 
\begin{equation}\label{Mnj-decomposition}
M_{n,j}=\begin{pmatrix}z&0\\ 0&1\end{pmatrix}\left(
I+\begin{pmatrix}0& r_{n,j}\\  r_{n,j}&0\end{pmatrix}\right)
\end{equation}
and, in general, 
\begin{equation}\label{noncommutativity}
\left(I+\begin{pmatrix}0&p\\ q&0\end{pmatrix}\right)\begin{pmatrix}s&0\\ 0&1\end{pmatrix}
=\begin{pmatrix}s&0\\ 0&1\end{pmatrix}\left(I+\begin{pmatrix}0&s^{-1}p\\ sq&0\end{pmatrix}\right).
\end{equation}
Some careful bookkeeping based on these observations yields formulas for the matrix entries of (\ref{OPUC-notation}). 
 In particular, the resulting formula for $\Psi_{n,n}^\ast$ is as follows. 
\begin{prop}\label{prop-psi-formula}
For each $1\leq j\leq n$, define $\kappa(s_1,\ldots,s_j)$ as in (\ref{kappa}), 
and denote
\[
R_{n,j}(\sigma)=\sum\limits_{1\leq\nu_1<\cdots<\nu_j\leq n}e^{2i\sigma\kappa(y_{n,\nu_1},\ldots,y_{n,\nu_j})}\left(\prod_{k=1}^j\alpha(y_{n,\nu_k})\right)\Delta_n^j.
\]
Then 
\[
\Psi_{n,n}^\ast\bigl(e^{2i\Delta_n\sigma}\bigr)=1+\sum_{j=1}^nR_{n,j}(\sigma). 
\]
\end{prop}

A key implication of (\ref{OPUC-1}) is that 
\[
\Psi_{n,n}^\ast\left(e^{2i\Delta_n\sigma}\right)\cong E^{(x_0,x)}_\alpha(\sigma).
\]
Based on Proposition~\ref{prop-psi-formula}, we prove a more precise result, as follows. 
\begin{lem}\label{lem-approximation}
Suppose $\alpha\in L^1_\real(X)$ is uniformly continuous, and let $K\subset \complex$ be an arbitrary compact set. Then, for every $x\in (x_0,x_1]$, 
\[
\Psi^\ast_{n,n}\left(e^{2i\Delta_n\sigma}\right)\rightarrow E^{(x_0,x)}_\alpha(\sigma)\quad\mbox{ as }\quad n\rightarrow\infty
\]
uniformly with respect to $\sigma\in K$. 
\end{lem}
\begin{pf}
Let $k_{j,\sigma}$ denote the integrand in $C_j$, as defined in (\ref{Cj}), so that 
\begin{equation}\label{kj}
k_{j,\sigma}(s_1,\ldots,s_j)=e^{2i\sigma\kappa(s_1,\ldots,s_j)}\displaystyle\prod_{\nu=1}^j\alpha(s_\nu).
\end{equation}
Observe that $R_{n,j}$ is a Riemann sum approximation to the integral $C_j$. 
Uniform continuity of $\alpha$ implies that for each fixed $j\geq 1$ the family of functions $k_{j,\sigma}$ $(\sigma\in K)$ is uniformly equicontinuous on $(x_0,x_1)^j$ (and in particular on the simplex $\mathscr{S}^{(x_0,x)}_{j,\rm{id}}$).  Therefore 
\[
R_{n,j}(\sigma)\rightarrow C_j(\sigma)\quad\mbox{ as }\quad n\rightarrow\infty
\]
uniformly with respect to $\sigma\in K$.  

Referring to (\ref{beta-sigma}) and (\ref{Cj-bound}), set $B=\max_{\sigma\in K}\beta_\sigma$ and $\alpha_{\max}=\sup_{x_0<x<x_1}|\alpha(x)|$.  Then 
\[
|C_j(\sigma)|\leq B\|\alpha\|^j/j!\leq \frac{B}{j!}\bigl((x_1-x_0)\alpha_{\max}\bigr)^j
\]
on $K$. Also,
\begin{equation}\label{R-bound}
\begin{split}
|R_{n,j}(\sigma)|&\leq \binom{n}{j}B(\Delta_n\alpha_{\max})^j\\
&<\frac{B}{j!}\bigl((x_1-x_0)\alpha_{\max}\bigr)^j.
\end{split}
\end{equation}
Thus, given $\varepsilon>0$ there exists $N$ sufficiently large that 
\[
2\sum_{j=N+1}^\infty \frac{B}{j!}\bigl((x_1-x_0)\alpha_{\max}\bigr)^j<\varepsilon/2. 
\]
And there exists $M$ sufficiently large that 
\[
\sum_{j=1}^N\left|R_{n,j}(\sigma)-C_j(\sigma)\right|<\varepsilon/2
\]
for every $n\geq M$ and $\sigma\in K$. Hence, for every $n\geq M$ and $\sigma\in K$, 
\[
\left|\Psi_{n,n}^\ast\bigl(e^{2i\Delta_n\sigma}\bigr)-E_{\alpha}^{(x_0,x)}(\sigma)\right|<\varepsilon,
\]
completing the proof. 
\end{pf}

Lemma~\ref{lem-approximation} extends immediately to the other entries of (\ref{OPUC-1}), by (\ref{dual-phi}) along with the fact that replacing $\alpha$ with its negative corresponds to interchanging $\Psi_{n,n}$ and $\Phi_{n,n}$.  Thus:
\begin{cor}\label{cor-approximation} 
If $\alpha\in L^1_{\real}(X)$ is uniformly continuous, then as $n\rightarrow\infty$,
\[
\begin{pmatrix}
\Psi_{n,n}\bigl(e^{2i\Delta_n\sigma}\bigr)&\Psi^\ast_{n,n}\bigl(e^{2i\Delta_n\sigma}\bigr)\\ -\Phi_{n,n}\bigl(e^{2i\Delta_n\sigma}\bigr)&\Phi_{n,n}^\ast\bigl(e^{2i\Delta_n\sigma}\bigr)
\end{pmatrix}\quad\rightarrow\quad
\begin{pmatrix}
e^{2i(x-x_0)\sigma}\overline{E^{(x_0,x)}_\alpha(\bar\sigma)}&E^{(x_0,x)}_\alpha(\sigma)\\ 
-e^{2i(x-x_0)\sigma}\overline{E^{(x_0,x)}_{-\alpha}(\bar\sigma)}&E^{(x_0,x)}_{-\alpha}(\sigma)
\end{pmatrix}
\]
uniformly on compact sets with respect to $\sigma\in\complex$. 
\end{cor}
\begin{remark}\label{remark-renormalization}
The partition width $\Delta_n=(x-x_0)/n$ serves as a renormalization factor in the left-hand matrix of Corollary~\ref{cor-approximation}.  The function $\Psi_{n,n}^\ast(e^{2i\Delta_n\sigma})$, for example, is periodic in $\sigma$, with period $\pi/\Delta_n$, which tends to $\infty$ as $n$ increases.  Thus the value of the harmonic exponential operator is a renormalized continuum limit of orthogonal polynomials. 
\end{remark}

\subsection{Further properties of the harmonic exponential\label{sec-properties-2}}
\begin{prop}\label{prop-real-part}
For every $\alpha\in L^1_\real(X)$, $x\in(x_0,x_1]$ and $\sigma\in\real$,
\[
\Re E^{(x_0,x)}_\alpha(\sigma)\overline{E^{(x_0,x)}_{-\alpha}(\sigma)}=1. 
\]
\end{prop}
\begin{pf}
Suppose first that $\alpha$ is uniformly continuous on $(x_0,x_1)$, and set 
\[
S_n(\sigma)=e^{-2i(x-x_0)\sigma}\det\begin{pmatrix}
\Psi_{n,n}\bigl(e^{2i\Delta_n\sigma}\bigr)&\Psi^\ast_{n,n}\bigl(e^{2i\Delta_n\sigma}\bigr)\\ -\Phi_{n,n}\bigl(e^{2i\Delta_n\sigma}\bigr)&\Phi_{n,n}^\ast\bigl(e^{2i\Delta_n\sigma}\bigr)
\end{pmatrix}.
\] 
If $\sigma\in\real$, then by Corollary~\ref{cor-approximation}
\begin{equation}\label{Sn-limit}
\begin{split}
S_n(\sigma)\rightarrow& e^{-2i(x-x_0)\sigma}\det\begin{pmatrix}
e^{2i(x-x_0)\sigma}\overline{E^{(x_0,x)}_\alpha(\sigma)}&E^{(x_0,x)}_\alpha(\sigma)\\ 
-e^{2i(x-x_0)\sigma}\overline{E^{(x_0,x)}_{-\alpha}(\sigma)}&E^{(x_0,x)}_{-\alpha}(\sigma)
\end{pmatrix}\\
&=2\Re E_\alpha^{(x_0,x)}(\sigma)\overline{E_{-\alpha}^{(x_0,x)}(\sigma)}
\end{split}
\end{equation}
uniformly on compact subsets of $\real$. Now, according to (\ref{Mnj}), 
\[
\det M_{n,j}=z\left(1-\bigl(\Delta_n\alpha(y_{n,j})\bigr)^2\right).
\]
Taking determinants in the formula (\ref{OPUC-notation}), with $j=n$, yields
\[
S_n(\sigma)=2\prod_{j=1}^n\left(1-\bigl(\Delta_n\alpha(y_{n,j})\bigr)^2\right). 
\]
Since $\alpha$ is bounded, $\bigl(\Delta_n\alpha(y_{n,j})\bigr)^2=O(n^{-2})$. Therefore,
\[
S_n(\sigma)\rightarrow 2\quad\mbox{ as }\quad n\rightarrow\infty.  
\]
The conclusion of the proposition then follows by (\ref{Sn-limit}). 

Next consider $\sigma\in\real$ and an arbitrary $\alpha\in L^1_{\real}(X)$.  Uniformly continuous functions are dense in $L^1_{\real}(X)$, so there exists a sequence $\alpha_n$ of uniformly continuous functions that converges to $\alpha$.  Proposition~\ref{prop-continuity}\ref{part-complex-continuity} then implies 
\[
1=\Re E^{(x_0,x)}_{\alpha_n}(\sigma)\overline{E^{(x_0,x)}_{-\alpha_n}(\sigma)}\quad\rightarrow\quad
\Re E^{(x_0,x)}_\alpha(\sigma)\overline{E^{(x_0,x)}_{-\alpha}(\sigma)}
\]
uniformly on compact subsets of $\real$, as $n\rightarrow\infty$. Hence $\Re E^{(x_0,x)}_\alpha(\sigma)\overline{E^{(x_0,x)}_{-\alpha}(\sigma)}=1$. 
\end{pf}
\begin{prop}\label{prop-zero-free} For every $\alpha\in L^1_{\real}(X)$, $x\in (x_0,x_1]$ and $\sigma\in\complex$ such that $\Im\sigma\geq 0$, 
each of $E^{(x_0,x)}_\alpha(\sigma)$, $E^{(x_0,x)}_{-\alpha}(\sigma)$ and $E^{(x_0,x)}_\alpha(\sigma)+E^{(x_0,x)}_{-\alpha}(\sigma)$ is nonzero. 
\end{prop}
\begin{pf}
A zero $\sigma$ of $E^{(x_0,x)}_\alpha$ or $E^{(x_0,x)}_{-\alpha}$ with $\Im\sigma=0$ would contradict Proposition~\ref{prop-real-part}. Proposition~\ref{prop-real-part} also implies that if $\Im\sigma=0$, then
\[
\bigl|E^{(x_0,x)}_\alpha(\sigma)+E^{(x_0,x)}_{-\alpha}(\sigma)\bigr|^2=\bigl|E^{(x_0,x)}_\alpha(\sigma)\bigr|^2+2\Re E^{(x_0,x)}_\alpha(\sigma)\overline{E^{(x_0,x)}_{-\alpha}(\sigma)}
+\bigl|E^{(x_0,x)}_{-\alpha}(\sigma)\bigr|^2>2,
\]
precluding any real zeros of $E^{(x_0,x)}_\alpha+E^{(x_0,x)}_{-\alpha}$.

Suppose next that $\alpha$ is uniformly continuous.  Then, on the upper half plane $\Im\sigma>0$, $E^{(x_0,x)}_\alpha(\sigma)$ is the uniform limit on compact sets of analytic functions $\Psi_{n,n}^\ast(e^{2i\Delta_n\sigma})$, by Corollary~\ref{cor-approximation}. If $\Im\sigma>0$, then $e^{2i\Delta_n\sigma}\in\ddd$, so these functions are zero free on the upper half plane by Proposition~\ref{prop-zero-free}. Hurwitz's theorem therefore implies $E^{(x_0,x)}_\alpha$ is either identically zero, or zero free, on the upper half plane. The former possibility is ruled out by the fact that $E^{(x_0,x)}_\alpha$ has no real zeros.  Exactly the same argument yields that $E^{(x_0,x)}_{-\alpha}$ and $E^{(x_0,x)}_\alpha+E^{(x_0,x)}_{-\alpha}$ are zero free on the upper half plane. 

Finally, if $\alpha\in L^1_{\real}(X)$ is arbitrary, then approximation in $L^1_{\real}(X)$ by uniformly continuous functions $\alpha_n$ yields the desired result, using Proposition~\ref{prop-continuity}\ref{part-complex-continuity} and Hurwitz's theorem. 
\end{pf}

\begin{lem}\label{lem-lower-bound}
Fix $\alpha\in L^1_{\real}(X)$, $x\in(x_0,x_1]$ and $\varepsilon>0$.  If $\alpha$ is uniformly continuous, then for every sufficiently large $n$, the inequality 
\[
\bigl|\Psi_{n,n}^\ast(z)+\Phi_{n,n}^\ast(z)\bigr|>\sqrt{2}-\varepsilon
\]
holds for every $z\in\overline{\ddd}$. 
\end{lem}
\begin{pf}
On one hand, if $|z|=1$, then (\ref{dual-phi}) implies
\[
\begin{split}
\det\begin{pmatrix}
\Psi_{n,n}(z)&\Psi^\ast_{n,n}(z)\\ -\Phi_{n,n}(z)&\Phi_{n,n}^\ast(z)
\end{pmatrix}&=\Psi_{n,n}(z)\Phi_{n,n}^\ast(z)+\Phi_{n,n}(z)\Psi^\ast_{n,n}(z)\\
&=z^n\left(\overline{\Psi^\ast_{n,n}(z)}\Phi_{n,n}^\ast(z)+\overline{\Phi_{n,n}^\ast(z)}\Psi^\ast_{n,n}(z)\right)\\
&=z^n2\Re\left(\Psi^\ast_{n,n}(z)\overline{\Phi_{n,n}^\ast(z)}\right).
\end{split}
\]
On the other hand, the right-hand side of (\ref{OPUC-notation}) (in the case $j=n$) yields
\[
\det\begin{pmatrix}
\Psi_{n,n}(z)&\Psi^\ast_{n,n}(z)\\ -\Phi_{n,n}(z)&\Phi_{n,n}^\ast(z)
\end{pmatrix}=z^n2\prod_{j=1}^n\left(1-\Delta_n^2\alpha_(y_{n,j})^2\right).
\]
Therefore
\[
\Re\left(\Psi^\ast_{n,n}(z)\overline{\Phi_{n,n}^\ast(z)}\right)=\prod_{j=1}^n\left(1-\Delta_n^2\alpha_(y_{n,j})^2\right),
\]
from which it follows that 
\begin{equation}\label{psi-phi-lower-bound}
\begin{split}
\left|\Psi_{n,n}^\ast(z)+\Phi_{n,n}^\ast(z)\right|^2&=
\bigl|\Psi_{n,n}^\ast(z)\bigr|^2+2\Re\left(\Psi^\ast_{n,n}(z)\overline{\Phi_{n,n}^\ast(z)}\right)
+\bigl|\Phi_{n,n}^\ast(z)\bigr|^2\\
&>2\prod_{j=1}^n\left(1-\Delta_n^2\alpha_(y_{n,j})^2\right)
\end{split}
\end{equation}
on the unit circle $|z|=1$. Now, $\Psi_{n,n}^\ast+\Phi_{n,n}^\ast$ is zero free on the $\overline{\ddd}$, by Proposition~\ref{prop-polynomials-zero-free}. The minimum modulus principle therefore implies (\ref{psi-phi-lower-bound}) holds for all $z\in\overline{\ddd}$.  Boundedness of $\alpha$ implies 
\[
\prod_{j=1}^n\left(1-\Delta_n^2\alpha_(y_{n,j})^2\right)\rightarrow1\quad\mbox{ as }\quad n\rightarrow\infty. 
\]
Hence, for any $\varepsilon>0$, the inequality 
\[
\sqrt{2}\prod_{j=1}^n\sqrt{1-\Delta_n^2\alpha_(y_{n,j})^2}>\sqrt{2}-\varepsilon
\]
holds for all sufficiently large $n$.  The desired conclusion then follows from (\ref{psi-phi-lower-bound}). 
\end{pf}
\begin{prop}\label{prop-lower-bound}
For every $\alpha\in L^1_{\real}(X)$, $x\in(x_0,x_1]$, and $\sigma\in\complex$ such that $\Im\sigma\geq0$, 
\[
\bigl|E^{(x_0,x)}_\alpha(\sigma)+E^{(x_0,x)}_{-\alpha}(\sigma)\bigr|\geq\sqrt{2}.
\]
\end{prop}
\begin{pf}
Suppose first that $\alpha$ is uniformly continuous, and hence bounded. If $\Im\sigma\geq0$, then for every $n\geq 1$, $z=e^{2i\Delta_n\sigma}\in\overline{\ddd}$.  Corollary~\ref{cor-approximation} and Lemma~\ref{lem-lower-bound} therefore imply that for any $\varepsilon>0$, 
\[
\sqrt{2}-\varepsilon\leq\lim_{n\rightarrow\infty}\bigl|\Psi_{n,n}^\ast\bigl(e^{2i\Delta_n\sigma}\bigr)+\Phi_{n,n}^\ast\bigl(e^{2i\Delta_n\sigma}\bigr)\bigr|=
\bigl|E^{(x_0,x)}_\alpha(\sigma)+E^{(x_0,x)}_{-\alpha}(\sigma)\bigr|.
\]
This yields the desired conclusion in the case where $\alpha$ is uniformly continuous. For arbitrary $\alpha\in L^1_\real(X)$, choose a sequence $\alpha_n$ of uniformly continuous functions converging to $\alpha$ in $L^1_\real(X)$.  If $\Im\sigma\geq0$, Proposition~\ref{prop-continuity}\ref{part-complex-continuity} implies 
\[
\sqrt{2}\leq\bigl|E^{(x_0,x)}_{\alpha_n}(\sigma)+E^{(x_0,x)}_{-\alpha_n}(\sigma)\bigr|
\rightarrow
\bigl|E^{(x_0,x)}_\alpha(\sigma)+E^{(x_0,x)}_{-\alpha}(\sigma)\bigr|\quad\mbox{ as }\quad n\rightarrow\infty,
\]
completing the proof. 
\end{pf}
\begin{prop}\label{prop-sech-i}
For every $\alpha\in L^1_\real(X)$ and $x\in(x_0,x_1]$, $\Sech^{(x_0,x)}_\alpha(i)>0$.  
\end{prop}
\begin{pf}
Evidently $\Sech^{(x_0,x)}_\alpha(i)\neq0$, since $\Cosh^{(x_0,x)}_\alpha$ is entire.  What needs to be shown is that $\Sech^{(x_0,x)}_\alpha(i)$ is real and positive. 
Note by Lemma~\ref{lem-approximation} that, if $\alpha$ is uniformly continuous,
\[
\begin{split}
E^{(x_0,x)}_\alpha(i)&=\lim_{n\rightarrow\infty}\Psi_{n,n}^\ast(e^{-2\Delta_n})\\
&=\lim_{n\rightarrow\infty}\Psi_{n,n}^\ast(1)\\
&=\lim_{n\rightarrow\infty}\prod_{j=1}^n\bigl(1+r_{n,j}\bigr),
\end{split}
\]
the latter equality by Proposition~\ref{prop-psi-formula}.  Thus $E^{(x_0,x)}_\alpha(i)\geq0$, since each factor $1+r_{n,j}>0$ for $n$ sufficiently large. Proposition~\ref{prop-zero-free} then implies $E^{(x_0,x)}_\alpha(i)>0$.  Replacing $\alpha$ with $-\alpha$ yields 
\[
E^{(x_0,x)}_{-\alpha}(i)=\lim_{n\rightarrow\infty}\prod_{j=1}^n\bigl(1-r_{n,j}\bigr), 
\]
from which it follows that $E^{(x_0,x)}_{-\alpha}(i)>0$ by a similar argument.  
Hence 
\[
\Sech^{(x_0,x)}_\alpha(i)=\frac{2}{E^{(x_0,x)}_\alpha(i)+E^{(x_0,x)}_{-\alpha}(i)}>0
\]
if $\alpha$ is uniformly continuous.  Approximation of arbitrary $\alpha\in L^1_\real(X)$ by a sequence $\alpha_n$ of uniformly continuous functions yields $\Sech^{(x_0,x)}_\alpha(i)>0$ by Proposition~\ref{prop-continuity}\ref{part-complex-continuity}, completing the proof.
\end{pf}

\section{Concluding remarks\label{sec-remarks}}

\subsection{Regularity of the potential $q$.\label{sec-regularity}}  Part of the significance of Theorem~\ref{thm-sech-outer} is to extend previously known characteristics of the scattering matrix \cite[Thm.~2, p.146]{DeTr:1979} to a setting of lower regularity.  Local integrability of $q$ required in \cite{DeTr:1979} corresponds to $\alpha$ being absolutely continuous, in contrast to Theorem~\ref{thm-sech-outer}, where $\alpha\in L^1_\real(X)$.  And of course in Theorem~\ref{thm-jumps}, where $\alpha$ can have Dirac delta singularities, corresponding to jumps in $\zeta$, the regularity is lower still.  Theorem~\ref{thm-sech-outer} also illustrates the utility of the harmonic exponential, which through its approximation by orthogonal polynomials allows one easily to deduce that $\Sech^{(x_0,x)}_\alpha(i)>0$.

\subsection{A continuum analogue of orthogonal polynomials\label{sec-continuum-analogue}} The analogy between the harmonic exponential and orthogonal polynomials can be fleshed out in heuristic terms as follows. The case $j+1=n$ of the classical recurrence (\ref{recurrence-psi-ast}), in which $y_{n,n}=x$ and $r_{n,n}=\Delta_n\alpha(x)$, can be rearranged in the form
\begin{equation}\label{modified-recurrence}
\frac{\Psi^\ast_{n,n}-\Psi_{n,n-1}^\ast}{e^{2i\Delta_n\sigma}\Delta_n}=\alpha(x)\Psi_{n,n-1},
\end{equation}
where $e^{2i\Delta_n\sigma}\cong1$ if $n$ is large relative to $\sigma$. 
In comparison, the harmonic exponential equation~\ref{x-derivative-dual}, with $x$ in place of $x_1$, 
\begin{equation}\label{modified-DE}
\partial_xE^{(x_0,x)}_\alpha=\alpha(x)\bigl(E_{\alpha}^{(x_0,x)}\bigr)^{\!\ast},
\end{equation}
can be seen as a continuous analogue of (\ref{modified-recurrence}), especially when one notes that 
\[
E^{(x_0,x)}_\alpha=\bigl(\bigl(E_{\alpha}^{(x_0,x)}\bigr)^{\!\ast}\bigr)^{\!\ast}. 
\]
In particular, the continuous analogue of the index $n$ is the spatial variable $x$, and the continuous 
analogue of the recursion coefficient $r_{n,n}=\Delta_n\alpha(y_{n,n})$ is the function $\alpha(x)$.  

The sequence $\Psi_{n,0},\ldots,\Psi_{n,n}$ is orthogonal on the unit circle with respect to the probability measure 
\begin{equation}\label{nun}
d\nu_n\bigl(e^{i\theta}\bigr)=\frac{\prod_{j=1}^n(1-| r_{n,j}|^2)}{2\pi\left|\Psi_{n,n}(e^{i\theta})\right|^2}\,d\theta\qquad (-\pi<\theta\leq \pi).
\end{equation}
Setting $\theta=2\Delta_n\sigma$ and $L_n=\pi/(2\Delta_n)$, the corresponding inner product, expressed in terms of $\sigma$, is 
\[
\langle f,g\rangle_{d\nu_n}=\frac{1}{2L_n}\int\limits_{-L_n}^{L_n}f(e^{2i\Delta_n\sigma})\overline{g(e^{2i\Delta_n\sigma})}\,\frac{\prod_{j=1}^n(1-| r_{n,j}|^2)}{\left|\Psi_{n,n}(e^{2i\Delta_n\sigma})\right|^2}\,d\sigma.
\]
Observe that the form
\[
\lim_{n\rightarrow\infty}\frac{1}{2L_n}\int\limits_{-L_n}^{L_n}\cdot\rule{5pt}{0pt} d\sigma
\]
corresponds precisely to the almost periodic scalar product of Besicovitch \cite{Be:1955},
\begin{equation}\label{ap}
\langle F,G\rangle_{\mbox{ap}}=\lim_{L\rightarrow\infty}\frac{1}{2L}\int\limits_{-L}^{L}F(\sigma)\overline{G(\sigma)}d\sigma,
\end{equation}
while the weight function in $d\nu_n$ stabilizes in the large $n$ limit as 
\[
\frac{\prod_{j=1}^n(1-| r_{n,j}|^2)}{\left|\Psi_{n,n}(e^{2i\Delta_n\sigma})\right|^2}\quad\rightarrow\quad 
\frac{1}{\bigl|\bigl(E^{(x_0,x)}_\alpha\bigr)^\ast(\sigma)\bigr|^2}=\frac{1}{\bigl|E^{(x_0,x)}_\alpha(\sigma)\bigr|^2}\qquad (\sigma\in\real).
\]
A consequence of the mean integral in (\ref{ap}) is that $\|F\|_{\mbox{ap}}=0$ if $F\in L^2(\real)$. Proposition~\ref{prop-hexp}\ref{hexp-limit} implies
\[
\frac{1}{\bigl|E^{(x_0,x)}_\alpha(\sigma)\bigr|^2}\quad\rightarrow\quad 1\quad\mbox{ as }\quad |\sigma|\rightarrow\infty\qquad(\sigma\in\real).
\] 
Assume $1/\bigl|E^{(x_0,x)}_\alpha(\sigma)\bigr|$ is an $L^2(\real)$-perturbation of the constant function $1$ (which is true if $\alpha\in L^2(X)$ \cite[Prop.~3.8(v)]{Gi:Ar2021}).  Then, from the point of view of almost periodic functions, the limiting weight is equivalent to the constant function 1.  Thus, orthogonality of two polynomials in the sequence $\Psi_{n,0},\ldots,\Psi_{n,n}$ with respect to $d\nu_n$ corresponds to orthogonality of two functions $\bigl(E^{(x_0,y_1)}_\alpha\bigr)^\ast$, $\bigl(E^{(x_0,y_2)}_\alpha\bigr)^\ast$ with respect to $\langle\cdot,\cdot\rangle_{\mbox{ap}}$, for some $x_0\leq y_1<y_2\leq x_1$. Since for $\sigma\in\real$, 
\[
\bigl(E^{(x_0,y_1)}_\alpha(\sigma)\bigr)^\ast\overline{\bigl(E^{(x_0,y_2)}_\alpha(\sigma)\bigr)^\ast}=e^{2iy_1\sigma}e^{-2iy_2\sigma}\overline{E^{(x_0,y_1}_\alpha(\sigma)}E^{(x_0,y_2)}_\alpha(\sigma)\quad\stackrel{\mbox{ap}}{\sim}\quad e^{2iy_1\sigma}e^{-2iy_2\sigma},
\]
this reduces to classic orthogonality of exponentials in the (non-separable!)~Besicovitch space of almost periodic functions.  Thus, from the perspective of the harmonic exponential, the continuum analogue of orthogonality on the circle is orthogonality on the line in the mean, as in (\ref{ap}).  A curious feature of this analogy is that the detailed structure of the weight function $1/\bigl|E^{(x_0,x)}_\alpha\bigr|^2$ is irrelevant to orthogonality. But of course it is essential in other respects, determining the scattering matrix $S$, for instance, as in Theorem~\ref{thm-schrodinger}.

From a more general analytic point of view, orthogonal polynomials on the unit circle elucidate the relationship of the Taylor coefficients of a bounded holomorphic function on the unit disk to its continued fraction expansion, as explained in \cite{SiOPUC1:2005,Kh:2008}. The harmonic exponential plays an analogous role in relating the scattering matrix of a Schr\"odinger operator in to its potential $q=\alpha^2-\alpha^\prime$. 

\subsection{Role of the Miura map\label{sec-miura}} In the present paper the Miura map carries the coefficient $\alpha$ in (\ref{wave-system}) to the potential $q$ in the classical time-independent eigenvalue form (\ref{schrodinger}) of the Schr\"odinger equation.  In its original formulation, however, the Miura map carries solutions of the modified Korteweg-de Vries equation (mKdV) to solutions of the original KdV equation; see \cite[\S2]{Mi:1968}.  At first blush, these are completely different situations.  But, of course, Schr\"odinger equations and the KdV hierarchy are famously connected via the Lax pair formulation, and the inverse scattering method for solution of the KdV hierarchy.  See \cite{AbKaNeSe:1974,BeDeTo:1988}.  In the simplest case, the classical Schr\"odinger equation corresponds to KdV, as follows. If one views a solution to KdV as a time-varying potential in (\ref{schrodinger}), then the associated time-varying reflection coefficient evolves according to a simple, linear evolution equation.  Solution of the latter, combined with inverse scattering, yields a solution to KdV. 
The system (\ref{wave-system}) corresponds to mKdV in a similar way.  Thus, although KdV does not enter into the present paper, the Lax pair formulation of the KdV hierarchy shows that the Miura map as it arises in the present article reflects the connection between mKdV and KdV as integrable systems.\footnote{Thanks to an anonymous referee for drawing attention to this point.}

\subsection{The notion of explicit solution\label{sec-explicit}}  The term ``explicit solution" is not precise, and gets used in various ways. For example, in \cite[Thm.~1]{GeGr:2015}, we see the inverse of a precisely defined operator referred to as an explicit solution to the cubic Szeg\H{o} equation.  In the present article the formula (\ref{hexp-formula}) describing the harmonic exponential is explicit in a stronger sense: not only can one analyze the formula directly, as is done in \S\ref{sec-analysis}, but truncation of the formula provides a direct means for numerical evaluation.  Indeed, the natural-seeming harmonic exponential equation (\ref{harmonic-exponential-equation}) suggests its solution is a basic object.  The series formulation (\ref{harmonic-exponential-equation}) is loosely akin to the scalar exponential series; a more explicit form is not to be expected.  Theorem~\ref{thm-general-solution} shows that the hyperbolic operators (\ref{hyperbolic-operators}) associated to the harmonic exponential yield an explicit representation for the solution to the Schr\"odinger equation naturally suited to the underlying algebraic structure of SU(1,1).

\appendix

\section{Proof of Proposition~\ref{prop-uniqueness}\label{sec-appendix-uniqueness}}

Fix $X=(x_0,x_1)$ and $\sigma\neq 0$, and suppose first that $\zeta$ is absolutely continuous on $X$, setting $\alpha=-\frac{1}{2}(\log\zeta)^\prime$.  Then any function $A$ that satisfies (\ref{wave-system}) is bounded on $X$ (in the sense that $a(x,\sigma)$ and $b(x,\sigma)$ are bounded).  This follows from the equivalence of (\ref{wave-system}) with (\ref{wave-fourier}), and the fact that if $\zeta$ is absolutely continuous, then any solution $u$ to (\ref{wave-fourier}) has the property that both $u$ and $u^\prime$ are absolutely continuous. Hence $a$ and $b$ are absolutely continuous and therefore bounded.  

Next, for any interval $Y=(y_0,y_1)$ such that $x_0\leq y_0<y_1\leq x_1$, let $V_Y$ denote the Banach space consisting of all bounded, measurable functions 
\begin{equation}\label{VY}
f=\binom{f_1}{f_2}:Y\rightarrow\complex^2\quad\mbox{ with norm }\quad \|f\|_Y=\sup_{x\in Y}\max\{|f_1(x)|,|f_2(x)|\}. 
\end{equation}
Consider the operator 
\begin{equation}\label{L}
L_{Y,\sigma}:V_Y\rightarrow V_Y,\qquad(L_{Y,\sigma}f)(x)=\int_{y_0}^x\frak{m}(s,\sigma)f(s)\,ds\quad (x\in Y),
\end{equation}
where
\[
\frak{m}(x,\sigma)=\begin{pmatrix}i\sigma&\alpha(x)\\ \alpha(x)&-i\sigma\end{pmatrix}. 
\]
Observe that 
\begin{equation}\label{L-bound}
\|L_{Y,\sigma}f\|_Y\leq \left((y_1-y_0)|\sigma|+\int_{y_0}^{y_1}|\alpha|\right)\|f\|_Y.  
\end{equation}
Therefore the operator norm on $L_{Y,\sigma}$ induced by $\|\cdot\|_Y$ is bounded by $$(y_1-y_0)|\sigma|+\int_{y_0}^{y_1}|\alpha|.$$  Now, since $\int_{x_0}^x|\alpha|$ is absolutely continuous for $x\in X$, for every $\epsilon>0$ there is a corresponding $\delta>0$ such that 
\begin{equation}\label{epsilon-delta}
y_1-y_0<\delta\quad\Rightarrow\quad (y_1-y_0)|\sigma|+\int_{y_0}^{y_1}|\alpha|<\epsilon. 
\end{equation}
Taking $\epsilon =1$, this implies that if the length of $Y\subset X$ is less than $\delta$, then the operator $L_{Y,\sigma}$ has norm less than 1. In this case the operator 
\begin{equation}\label{invertible}
1-L_{Y,\sigma}:V_Y\rightarrow V_Y
\end{equation}
is invertible, and $f=L_{Y,\sigma}f$ only if $f=0$.  It follows that $f=L_{X,\sigma}f$ only if $f=0$.  To see this, cover $X$ with a sequence of small open intervals $Y_j$ $(1\leq j\leq n)$, indexed and arranged such that the sequences of left and right endpoints are both increasing.  In particular, the left end point of $Y_1$ has to be $x_0$, and the right endpoint of $Y_n$ has to be $x_1$.  If $f=L_{X,\sigma}f$, then necessarily $f|_{Y_1}=L_{Y_1,\sigma}(f|_{Y_1})$ and so $f|_{Y_1}=0$.  Given that $f|_{Y_1}=0$, the assumption $f=L_{X,\sigma}f$ implies in turn that $f|_{Y_2}=L_{Y_2,\sigma}(f|_{Y_2})$, since, for $x\in Y_2$, 
\[
L_{X,\sigma}f(x)=L_{Y_2,\sigma}f(x). 
\]
Since $Y_2$ is small enough, $f|_{Y_2}=0$ also.  Continuing in this way exhausts all of $X$, proving $f|_{Y_j}=0$ for each $j$, whence $f=0$.  

The system (\ref{wave-system}) is equivalent to the integral equation
\begin{equation}\label{integral-form}
A(x,\sigma)=A(x_0+,\sigma)-\int_{x_0}^x\frak{m}(\xi,\sigma)A(\xi,\sigma)\,d\xi. 
\end{equation}
Suppose $A$ satisfies (\ref{integral-form}) with $A(x_0+,\sigma)=0$.  Then $A=L_{X,\sigma}A$, and so $A=0$ by the foregoing argument.  Linearity then implies equation (\ref{integral-form}) with arbitrary $A(x_0+,\sigma)$ has at most one solution.  The equivalent differential equation (\ref{wave-system}) therefore also has at most one solution, as does the original equation (\ref{wave-fourier}), assuming $u(x_0+,\sigma)$ and $u^\prime(x_0+,\sigma)$ are fixed.  

Suppose next that $\zeta$ is piecewise absolutely continuous.  Let $y_1,\ldots,y_n$ be the jump points of $\zeta$, labelled such that 
\[
x_0=y_0<y_1<\cdots y_n<y_{n+1}=x_1.  
\]
Then on any interval $Y_j=(y_{j-1},y_j)$ $(1\leq j\leq n+1)$, there is at most one solution to the restriction of (\ref{wave-system}) to $Y_j$ consistent with a given value of $A(y_{j-1}+,\sigma)$.  In particular the restriction of $A$ to $Y_1$ is determined by $A(x_0+,\sigma)$ and determines $A(y_1-,\sigma)$.  Then $A(y_1+,\sigma)$ is determined from $A(y_1-,\sigma)$ by the jump condition (\ref{scattering-at-x}).  This determines $A$ on $Y_2$.  Continuing this procedure determines $A$ on all of $X$.   

Lastly, in the case $\sigma=0$, the original equation (\ref{wave-fourier}) can be solved directly by integration, with the solution being determined by $u(x_0+,0)$ and $u^\prime(x_0+,0)$.

%%%%%%%%%%%%%%%%
%%%%%%%%%%%%%%%%
%\bibliographystyle{abbrv}
%\bibliography{References_24.7.2020}
%%%%%%%%%%%%%%%%
%%%%%%%%%%%%%%%%

%\bibitem{Gi:Ar2021}
%P.~C. Gibson.
%\newblock Scattering on the line via singular approximation.
%\newblock arXiv:2108.09799 [math.AP], 88 pp., 2021.
%

\end{document}